\documentclass[12pt]{amsart}
\usepackage{amsmath}
\usepackage{amssymb}
\usepackage{amsfonts}
\usepackage{amsthm}
\usepackage{amscd}
\usepackage{bm}
\newtheorem{thrm}{Theorem}
\newtheorem{coro}{Corollary}

\newtheorem{thm}{Theorem}[section]
\newtheorem{lem}[thm]{Lemma}
\newtheorem{pro}[thm]{Proposition}
\newtheorem{cor}[thm]{Corollary}

\newtheorem{conj}[thm]{Conjecture}

\newcommand{\Z}{\mathbf{Z}}
\newcommand{\Q}{\mathbf{Q}}

\begin{document}
\title{On the low dimensional cohomology groups of \\ the IA-automorphism group of a free group of \\ rank three}
\address{Department of Mathematics, Faculty of Science Division II, Tokyo University of Science,
         1-3, Kagurazaka, Shinjuku-ku, Tokyo 162-8601, Japan}
\email{takao@rs.tus.ac.jp}
\maketitle
\begin{center}
{\textsf{Dedicated to Professor Shigeyuki Morita on the occasion of his 70th birthday}}

\vspace{1em}

{\sc Takao Satoh} \\

\vspace{0.5em}

{\footnotesize Department of Mathematics, Faculty of Science Division II, Tokyo University of Science, \\
         1-3, Kagurazaka, Shinjuku-ku, Tokyo 162-8601, Japan}
\end{center}

\begin{abstract}
In this paper we study the structure of the rational cohomology groups of the IA-automorphism group $\mathrm{IA}_3$ of a free group of rank three
by using combinatorial group theory and representation theory. In particular, we detect non-trivial irreducible component in the
second cohomology group of $\mathrm{IA}_3$, which does not contained in the image of the cup product map of the first cohomology groups.
We also show that the image of the triple cup product map of the first cohomology groups in the third cohomology group is trivial.
As a corollary, we obtain that the fourth term of the lower central series of $\mathrm{IA}_3$ has finite index in that of
the Andreadakis-Johnson filtration.
\end{abstract}
\section{Introduction}\label{Int}

Let $F_n$ be a free group of rank $n \geq 2$ with basis $x_1, \ldots , x_n$, and $\mathrm{Aut}\,F_n$ the automorphism group of $F_n$.
As far as we know, the first contribution to the study of the (co)homology groups of $\mathrm{Aut}\,F_n$ was given by Nielsen \cite{Ni1} in 1924,
who showed $H_1(\mathrm{Aut}\,F_n,\Z)=\Z/2\Z$ for $n \geq 2$ by using a presentation for $\mathrm{Aut}\,F_n$.
Now we have a broad range of results for the (co)homology groups of $\mathrm{Aut}\,F_n$ with trivial coefficients by many authors.
In 1984, Gersten \cite{Ger} showed $H_2(\mathrm{Aut}\,F_n,\Z)=\Z/2\Z$ for $n \geq 5$
In 1980s, by introducing the Outer space, 
Culler and Vogtmann \cite{CV1} made a breakthrough in computation of homology groups of the outer automorphism groups
of free groups. To put it briefly, the Outer space is an analogue of the Teichm$\ddot{\mathrm{u}}$ller space on which the mapping class of a surface
naturally acts.
By using the geometry of the Outer space,
Hatcher and Vogtmann \cite{HV2} computed $H_4(\mathrm{Aut} \,F_4, \Q)= \Q$.
On the other hand, by using sophisticated homotopy theory,
Galatius \cite{Gal} showed that the stable integral homology groups of
$\mathrm{Aut}\,F_n$ are isomorphic to those of the
symmetric group $\mathfrak{S}_n$ of degree $n$.
In particular, the stable rational homology groups $H_q(\mathrm{Aut}\,F_n,\Q)$ are trivial for $n \geq 2q+1$.
This result is a quite contrast to the case of the mapping class groups of surfaces.
Intuitively, we can see this from the fact that the free group has no geometric extra structure unlike surface groups.

\vspace{0.5em}

With respect to unstable cohomology groups, $\mathrm{Aut}\,F_n$ behaves in much different and mysterious way.
The unstable cohomology groups of the (outer) automorphism groups of free groups have also been studied by many authors.
The Outer space is of course a powerful tool for computation of unstable cohomology groups. For example,
Brady \cite{Bra} computed the integral cohomology groups of $\mathrm{Out}\,F_3$ in 1993, 
Gerlitz showed $H_7(\mathrm{Aut}\,F_5,\Q)=\Q$ in 2002, and Ohashi \cite{Oha} computed $H_8(\mathrm{Out} \,F_6, \Q)=\Q$ in 2007.
On the other hand, 
in 1999, Morita \cite{Mo2} constructed a series of unstable homology classes of $\mathrm{Out}\,F_n$
with Kontsevich's results \cite{Kon1} and \cite{Kon2}. (See also \cite{Mo3}.)
These homology classes are called the Morita classes. It is known that the first and the second one are non-trivial,
and hence are generators of $H_4(\mathrm{Out} \,F_4, \Q)$ and $H_8(\mathrm{Out} \,F_6, \Q)$ respectively.
(See \cite{Mo3} and \cite{CoV} respectively.)
Furthermore, Morita-Sakasai-Suzuki \cite{MSS} showed that $\mathrm{Out}\,F_n$ has many non-trivial unstable homology classes with trivial coefficients.
In spite of these intense studies, it seems that the structure of unstable (co)homology classes of $\mathrm{Out}\,F_n$ is still complicated.
We should remark that in \cite{CHKV}, Conant-Hatcher-Kassabov-Vogtmann gave a construction of many non-trivial unstable homology classes of
$\mathrm{Aut}\,F_n$ and $\mathrm{Out}\,F_n$, and studies the Morita classes.

\vspace{0.5em}

Let $H$ be the abelianization of $F_n$, and $\mathrm{IA}_n$ the kernel of the natural homomorphism $\mathrm{Aut}\,F_n \rightarrow \mathrm{Aut}\,H$
induced from the abelianization homomorphism $F_n \rightarrow H$. The group $\mathrm{IA}_n$ is called the IA-automorphism group of $F_n$.
By observing the spectral sequence of the group extension of $\mathrm{IA}_n$ by $\mathrm{Aut}\,H$, 
we see that the cohomology groups of $\mathrm{IA}_n$
are closely related to those of $\mathrm{Aut}\,F_n$.
However, the structure of the cohomology groups of $\mathrm{IA}_n$ is far from well-understood in contrast to those of $\mathrm{Aut}\,F_n$.
To our best knowledge, in the (co)homology groups of $\mathrm{IA}_n$, completely determined and explicitly written down one is
only the first integral homology group $H_1(\mathrm{IA}_n, \Z)$,
which is obtained by Cohen-Pakianathan \cite{Co1, Co2}, Farb \cite{Far} and Kawazumi \cite{Kaw} independently.
(See Subsection {\rmfamily \ref{Ss-IA}} for details.)
Krsti\'{c} and McCool \cite{Krs} showed that $\mathrm{IA}_3$ is not finitely presentable. This shows that there is a possibility that the second homology group
$H_2(\mathrm{IA}_3,\Z)$ is not finitely generated. 
In fact, this follows by a work of Bestvina, Bux and Margalit \cite{Bes}. By using the Outer space,
they showed that the quotient group of $\mathrm{IA}_n$ by the inner automorphism group $\mathrm{Inn}\,F_n$ has a $2n-4$-dimensional
Eilenberg-Maclane space, and that $H_{2n-4}(\mathrm{IA}_n/\mathrm{Inn}\,F_n, \Z)$ is not finitely generated.
For $n \geq 4$, it is not known whether $\mathrm{IA}_n$ is finitely presentable or not. Namely, at the present stage,
even $H_2(\mathrm{IA}_n,\Z)$ is not determined explicitly.
Pettet \cite{Pet} determined the image of the rational cup product of the first cohomologies in $H^2(\mathrm{IA}_n,\Q)$,
and gave its irreducible $\mathrm{GL}$-decomposition.
Furthermore, Day and Putman \cite{D-P} obtained
an explicit finite set of generators for $H_2(\mathrm{IA}_n, \Z)$ as a $\mathrm{GL}(n, \Z)$-module.

\vspace{0.5em}

In this paper, we mainly study the second rational cohomology group $H^2(\mathrm{IA}_n,\Q)$ for the case where $n=3$.
In particular, we detect a new $\mathrm{GL}(3,\Q)$-irreducible module in
$H^2(\mathrm{IA}_3,\Q)$ by using combinatorial group theory and representation theory.
By Pettet \cite{Pet}, the $\mathrm{GL}(3,\Q)$-irreducible decomposition of the image of the cup product
$\cup_{\Q} : \Lambda^2 H^1(\mathrm{IA}_3,\Q) \rightarrow H^2(\mathrm{IA}_3,\Q)$ is obtained.
As a corollary to Theorem {\rmfamily \ref{T-main1}} below,
we obtain the following.
\begin{thrm}[$=$ Theorem {\rmfamily \ref{T-main1}}]\label{I-T-1}
The quotient module $H^2(\mathrm{IA}_3,\Q)/\mathrm{Im}(\cup_{\Q})$ contains
a $\mathrm{GL}(3,\Q)$-irreducible representation $[3,2,-2]^*=[2,-2,-3] \cong D^{-3} \otimes_{\Q} [5,1]$.
\end{thrm}
\noindent
(For notations, see Subsection {\rmfamily \ref{Ss-Rep}}.) 

\vspace{0.5em}

In order to show Theorem {\rmfamily \ref{I-T-1}}, we use our previous results about
the Andreadakis-Johnson filtration $\mathrm{IA}_n = \mathcal{A}_n(1) \supset \mathcal{A}_n(2) \supset \cdots$
and the Johnson homomorphisms of $\mathrm{Aut}\,F_3$.
Historically, the Andreadakis-Johnson filtration was originally introduced by Andreadakis \cite{And} in the 1960s.
In 1980s, Johnson used this type of filtration to study the group structure of the mapping class groups of surfaces.
Andreadakis conjectured that the filtration $\mathrm{IA}_n = \mathcal{A}_n(1) \supset \mathcal{A}_n(2) \supset \cdots$ coincides
with the lower central series $\mathrm{IA}_n = \mathcal{A}_n'(1) \supset \mathcal{A}_n'(2) \supset \cdots$.
Andreadakis showed that this conjecture is true for $n=2$ and any $k \geq 2$, and $n=3$ and $k \leq 3$.
Bachmuth \cite{Ba2} showed $\mathcal{A}_n'(2) = \mathcal{A}_n(2)$ for any $n \geq 2$.
This result is also induced from the fact that the first Johnson homomorphism is the abelianization of $\mathrm{IA}_n$ by independent works 
Cohen-Pakianathan \cite{Co1, Co2}, Farb \cite{Far} and Kawazumi \cite{Kaw}.
Pettet \cite{Pet} showed that $\mathcal{A}_n'(3)$ has at most finite index in $\mathcal{A}_n(3)$ for any $n \geq 4$.
Recently, in our previous paper \cite{S25}, we showed that $\mathcal{A}_n'(3) = \mathcal{A}_n(3)$ for any $n \geq 3$.
By using a computer, Bartholdi \cite{Bar} showed that Andreadakis conjecture is not true for 
the case where $n=3$ and degree $7$.
In general, the Andreadakis conjecture is still open problem in a stable range
\footnote{Bartholdi \cite{Bar} asserted that the \lq\lq rational" version of the Andreadakis conjecture is true.
Andrew Putman, however, pointed out some gaps in Bartholdi's argument. Then,
in their communications, Bartholdi agreed them, and he published an erratum \cite{Bar2}.
Thus the conjecture is still open now.}.
In this paper, from our computation in the proof of Theorem {\rmfamily \ref{I-T-1}}, as a corollary, we also obtain the following.
\begin{coro}[$=$ Corollary {\rmfamily \ref{C-And}}]
$\mathcal{A}_3(4)/\mathcal{A}_3'(4)$ is finite.
\end{coro}
\noindent
We also remark that the affirmative answer to the Andreadakis conjecture resricted to certain subgroups of $\mathrm{Aut}\,F_n$ were given by \cite{S07} and \cite{S05},
and these works are generalized systematically by recent notable works of Darn\'{e} \cite{Da1, Da2}.

\vspace{0.5em}
Finally, we consider the third rational cohomology group $H^3(\mathrm{IA}_3,\Q)$.
The results by work of Bestvina, Bux and Margalit \cite{Bes} as mentioned above, we see that $H^3(\mathrm{IA}_3,\Q)$ is infinitely generated.
The following theorem shows that non-trivial elements in $H^3(\mathrm{IA}_3,\Q)$ cannot be detected by the triple cup product
of the first cohomology group of $\mathrm{IA}_3$.
\begin{thrm}[$=$ Theorem {\rmfamily \ref{T-main2}}]\label{I-T-2}
The image of the triple cup product
\[ \cup_{\mathrm{IA}_3}^3 : \Lambda^3 H^1(\mathrm{IA}_3, \Q) \rightarrow H^3(\mathrm{IA}_3, \Q) \]
is trivial.
\end{thrm}

\vspace{0.5em}

We remark that the arguments and techniques which we use in this paper are applicable to study the cohomology groups of $\mathrm{IA}_n$
for general $n \geq 4$. However, the amount of calculation and the complexity vastly increase with the increasing $n$.
In the present paper, we give the first combinatorial group theoretic approach to the study of the low dimensional cohomology
groups of the IA-automorphism groups of free groups.

\tableofcontents

\section{Preliminaries}\label{S-Pre}

In this section, after fixing some notation and conventions,
we recall the IA-automorphism groups, Free Lie algebras,
the Andreadakis-Johnson filtration and the representation theory of the general linear group over $\Q$.

\subsection{Notation and conventions}\label{S-Not}
\hspace*{\fill}\ 

\vspace{0.5em}

Let $G$ be a group.

\begin{itemize}
\item The abelianization of $G$ is denoted by $G^{\mathrm{ab}}$.
\item The automorphism group $\mathrm{Aut}\,G$ acts on $G$ from the right. For any $\sigma \in \mathrm{Aut}\,G$ and $x \in G$,
      the action of $\sigma$ on $x$ is denoted by $x^{\sigma}$.
\item For a normal subgroup $N$, we often denote the coset class of an element $g \in G$ by the same $g$ in the quotient group $G/N$ if
      there is no confusion.
\item For elements $x$ and $y$ of $G$, the commutator bracket $[x,y]$ of $x$ and $y$
      is defined to be $[x,y]:=xyx^{-1}y^{-1}$. Then for any $x, y, z \in G$, we have
\begin{eqnarray}
  & & [xy,z] =[x,[y,z]][y,z][x,z], \hspace{2em} [x,yz] = [x,y][x,z][[z,x],y] \label{eq-akagi} \\
  & & [x^{-1},z]=[[x^{-1}, z],x][x,z]^{-1}, \hspace{1em} [x,y^{-1}] = [x,y]^{-1} [y,[y^{-1}, x]] \label{eq-kaga}
\end{eqnarray}
      For elements $g_1, \ldots, g_k \in G$, a simple $k$-fold commutator
      $[[ \cdots [[ y_{1},y_{2}],y_{3}], \cdots ], y_{k}]$ is denoted by $[y_{1}, y_{2}, \cdots, y_{k}]$ for simplicity.
      For subgroups $H$ and $K$ of $G$, we denote by $[H,K]$ the commutator subgroup of $G$
      generated by $[h, k]$ for $h \in H$ and $k \in K$.
\item For any $\Z$-module $M$, we denote $M \otimes_{\Z} \Q$ by the symbol obtained by attaching a subscript or a superscript $\Q$ to $M$,
      like $M_{\Q}$ or $M^{\Q}$. Similarly, for any $\Z$-linear map $f: A \rightarrow B$,
      the induced $\Q$-linear map $f \otimes \mathrm{id}_{\Q} : A_{\Q} \rightarrow B_{\Q}$ is denoted by $f_{\Q}$
      or $f^{\Q}$.
\end{itemize}

\subsection{IA-automorphism groups}\label{Ss-IA}
\hspace*{\fill}\ 

\vspace{0.5em}

Fix a basis $x_1, \ldots , x_n$ of a free group $F_n$.
We denote by $H$ the abelianization $H_1(F_n, \Z)$ of $F_n$. Let $\rho : \mathrm{Aut}\,F_n \rightarrow \mathrm{Aut}\,H$
be the natural homomorphism induced from
the abelianization of $F_n$. We identify $\mathrm{Aut}\,H$ with the general linear group $\mathrm{GL}(n,\Z)$ by
fixing the basis of $H$ induced from the basis $x_1, \ldots , x_n$ of $F_n$.
The kernel $\mathrm{IA}_n$ of $\rho$ is called the IA-automorphism group of $F_n$.
It is clear that the inner automorphism group $\mathrm{Inn}\,F_n$
of $F_n$ is contained in $\mathrm{IA}_n$. Nielsen \cite{Ni0} showed that $\mathrm{IA}_2 = \mathrm{Inn}\,F_2$.
For $n \geq 3$, $\mathrm{IA}_n$ is much larger than $\mathrm{Inn}\,F_n$.
In fact, Magnus \cite{Mag} showed that $\mathrm{IA}_n$ is finitely generated by automorphisms
\[ K_{ij} : x_t \mapsto \begin{cases}
               {x_j}^{-1} x_i x_j, & t=i, \\
               x_t,                & t \neq i
              \end{cases}\]
for distinct $i$, $j \in \{ 1, 2, \ldots , n \}$ and
\[  K_{ijl} : x_t \mapsto \begin{cases}
               x_i [x_j, x_l], & t=i, \\
               x_t,            & t \neq i
              \end{cases}\] 
for distinct $i$, $j$, $l \in \{ 1, 2, \ldots , n \}$ such that $j>l$. 
For any $1 \leq i \leq n$, set $\iota_i := K_{1i}K_{2i} \cdots K_{ni}$. Namely, $\iota_i$ is the inner automorphism of $F_n$ given by
$x \mapsto x_i^{-1} x x_i$.

\vspace{0.5em}

Cohen-Pakianathan \cite{Co1, Co2}, Farb \cite{Far} and Kawazumi \cite{Kaw} independently showed
\begin{equation}\label{CPFK}
H_1(\mathrm{IA}_n, \Z) \cong H^* \otimes_{\Z} \Lambda^2 H
\end{equation}
as a $\mathrm{GL}(n,\Z)$-module where $H^*:= \mathrm{Hom}_{\Z}(H,\Z)$ is the $\Z$-linear dual group of $H$.
This fact is obtained from the above result of Magnus and the first Johnson homomorphism which is defined below.

\subsection{Free Lie algebra generated by $H$}\label{Ss-FL}
\hspace*{\fill}\ 

\vspace{0.5em}

Let $\Gamma_n(1) \supset \Gamma_n(2) \supset \cdots$ be the lower central series of a free group $F_n$ defined by the rule
\[ \Gamma_n(1):= F_n, \hspace{1em} \Gamma_n(k) := [\Gamma_n(k-1),F_n], \hspace{1em} k \geq 2. \]
We denote by $\mathcal{L}_n(k) := \Gamma_n(k)/\Gamma_n(k+1)$ the graded quotient of the lower central series of $F_n$,
and by $\mathcal{L}_n := {\bigoplus}_{k \geq 1} \mathcal{L}_n(k)$ the associated graded sum. See $\mathcal{L}_n(1)=H$.
Since the group $\mathrm{Aut}\,F_n$ naturally acts on $\mathcal{L}_n(k)$ for each $k \geq 1$,
and since $\mathrm{IA}_n$ acts on it trivially,
the action of $\mathrm{GL}(n,\Z)$ on each $\mathcal{L}_n(k)$ is well-defined.
Furthermore, the graded sum $\mathcal{L}_n$ naturally has a graded Lie algebra structure induced from
the commutator bracket on $F_n$, and is isomorphic to the free Lie algebra generated by $H$.
(See \cite{MKS} and \cite{Reu} for basic material concerning the free Lie algebra.)
It is well known due to Witt \cite{Wit} that each $\mathcal{L}_n(k)$ is a $\mathrm{GL}(n,\Z)$-equivariant
free abelian group of rank
\begin{equation}\label{ex-witt}
 r_n(k) := \frac{1}{k} \sum_{d | k} \mu(d) n^{\frac{k}{d}}
\end{equation}
where $\mu$ is the M$\ddot{\mathrm{o}}$bius function, and $d$ runs over all positive divisors of $k$.
For example, the $\mathrm{GL}(n,\Z)$-module structure of $\mathcal{L}_n(k)$ for $1 \leq k \leq 3$ is given by
\[\begin{split}
  \mathcal{L}_n(1) & = H, \hspace{1em} \mathcal{L}_n(2) = \Lambda^2 H, \\
  \mathcal{L}_n(3) & = (H \otimes_{\Z} \Lambda^2 H) \big{/} \langle x \otimes y \wedge z + y \otimes z \wedge x
   + z \otimes x \wedge y \,\, | \,\, x,y,z \in H \rangle.
  \end{split}\]

\vspace{0.5em}

Hall \cite{Hal} constructed an explicit basis of $\mathcal{L}_n(k)$.
More precisely, he introduced basic commutators of $F_n$, and showed that (the coset classes of) basic commutators of weight $k$
form a basis of $\mathcal{L}_n(k)$. For example, basic commutators of weight less than four are listed below.
\vspace{0.5em}
\begin{center}
{\renewcommand{\arraystretch}{1.5}
\begin{tabular}{|c|l|l|} \hline
  $k$  & basic commutators  &     \\ \hline
  $1$  & $x_1, \ldots, x_n$ &     \\ \hline
  $2$  & $[x_i,x_j]$        & $j<i$    \\ \hline
  $3$  & $[x_i, x_j, x_l]$  & $i > j \leq l$    \\ \hline
\end{tabular}}
\end{center}
\vspace{0.5em}
(See also \cite{Ha2} for details for the basic commutators of the free groups.)

\subsection{Andreadakis-Johnson filtration}\label{Ss-A-J}
\hspace*{\fill}\ 

\vspace{0.5em}

In this subsection, we recall the Andreadakis-Johnson filtration and Johnson homomorphisms of $\mathrm{Aut}\,F_n$.
For each $k \geq 1$, the action of $\mathrm{Aut}\,F_n$ on the nilpotent quotient group $F_n/\Gamma_n(k+1)$ of $F_n$ induces a homomorphism
\[ \mathrm{Aut}\,F_n \rightarrow \mathrm{Aut}(F_n/\Gamma_n(k+1)). \]
We denote its kernel by $\mathcal{A}_n(k)$. Then the groups $\mathcal{A}_n(k)$ define a descending central filtration
\[ \mathrm{IA}_n = \mathcal{A}_n(1) \supset \mathcal{A}_n(2) \supset \mathcal{A}_n(3) \supset \cdots \]
of $\mathrm{IA}_n$. We call this filtration the Andreadakis-Johnson filtration of $\mathrm{Aut}\,F_n$.
Andreadakis showed that
\begin{thm}[Andreadakis \cite{And}]\label{T-And} \quad
\begin{enumerate}
\item For any $k$, $l \geq 1$, $\sigma \in \mathcal{A}_n(k)$ and $x \in \Gamma_n(l)$, $x^{-1} x^{\sigma} \in \Gamma_n(k+l)$.
\item For any $k$ and $l \geq 1$, $[\mathcal{A}_n(k), \mathcal{A}_n(l)] \subset \mathcal{A}_n(k+l)$.
\item $\displaystyle \bigcap_{k \geq 1} \mathcal{A}_n(k) =1$.
\end{enumerate}
\end{thm}
For each $k \geq 1$, the group $\mathrm{Aut}\,F_n$ acts on $\mathcal{A}_n(k)$ by conjugation, and
it naturally induces an action of $\mathrm{GL}(n,\Z)=\mathrm{Aut}\,F_n/\mathrm{IA}_n$ on
the graded quotients $\mathrm{gr}^k (\mathcal{A}_n) := \mathcal{A}_n(k)/\mathcal{A}_n(k+1)$ by Part (2) of Theorem {\rmfamily \ref{T-And}}.
The graded sum ${\mathrm{gr}}(\mathcal{A}_n) := \bigoplus_{k \geq 1} {\mathrm{gr}}^k (\mathcal{A}_n)$ has a
graded Lie algebra structure induced from the commutator bracket on $\mathrm{IA}_n$.

\vspace{0.5em}

For each $k \geq 1$, define the homomorphism
$\tilde{\tau}_k : \mathcal{A}_n(k) \rightarrow \mathrm{Hom}_{\Z}(H, {\mathcal{L}}_n(k+1))$ by
\[ \sigma \hspace{0.3em} \mapsto \hspace{0.3em} (x \pmod{\Gamma_n(2)} \mapsto x^{-1} x^{\sigma} \pmod{\Gamma_n(k+2)}), \hspace{1em} x \in F_n. \]
Then the kernel of $\tilde{\tau}_k$ is $\mathcal{A}_n(k+1)$. 
Hence it induces the injective homomorphism
\[ \tau_k : \mathrm{gr}^k (\mathcal{A}_n) \hookrightarrow \mathrm{Hom}_{\Z}(H, \mathcal{L}_n(k+1))
       = H^* \otimes_{\Z} \mathcal{L}_n(k+1). \]
The homomorphisms $\tilde{\tau}_k$ and ${\tau}_{k}$ are called the $k$-th Johnson homomorphisms of $\mathrm{Aut}\,F_n$.
Each ${\tau}_{k}$ is a $\mathrm{GL}(n,\Z)$-equivariant homomorphism.
For the Magnus generators of $\mathrm{IA}_n$, their images by $\tau_1$ are given by
\begin{equation}
  \tau_1(K_{ij}) = x_i^* \otimes [x_i,x_j], \hspace{1em} \tau_1(K_{ijl}) = x_i^* \otimes [x_j,x_l]. \label{eq-jm}
\end{equation}
Hence $\tau_1$ is surjective. From this fact, we see that the first Johnson homomorphism induces the abelianization of $\mathrm{IA}_n$.
(For details about the Johnson homomorphisms, see \cite{S06}, \cite{TS1} and \cite{TS2} for example.)

\vspace{0.5em}

Let $\mathrm{IA}_n = \mathcal{A}_n'(1) \supset \mathcal{A}_n'(2) \supset \cdots$ be the lower central series of $\mathrm{IA}_n$,
and set $\mathrm{gr}^k(\mathcal{A}_n'):= \mathcal{A}_n'(k)/\mathcal{A}_n'(k+1)$ for each $k$.
Since the Andreadakis-Johnson filtration is central by Part (2) of Theorem {\rmfamily \ref{T-And}},
we see $\mathcal{A}_n'(k) \subset \mathcal{A}_n(k)$ for any $k \geq 1$. Then we have the following conjecture.
\begin{conj}[Andreadakis's conjecture]
For any $n \gg k \geq 1$, $\mathcal{A}_n'(k) = \mathcal{A}_n(k)$.
\end{conj}

\vspace{0.5em}

The restriction of $\tilde{\tau}_k : \mathcal{A}_n(k) \rightarrow \mathrm{Hom}_{\Z}(H, \mathcal{L}_n(k+1))$
to $\mathcal{A}_n'(k)$ induces the homomorphism
\[ \tau_k' : \mathrm{gr}^k (\mathcal{A}_n') \rightarrow \mathrm{Hom}_{\Z}(H, \mathcal{L}_n(k+1)). \]
In this paper, by abuse of language, we also call $\tau_k'$ the $k$-th Johnson homomorphism of $\mathrm{Aut}\,F_n$.
We can see that each $\tau_k'$ is $\mathrm{GL}(n,\Z)$-equivariant by the same way as $\tau_k$.
Then we have the exact sequence
\[ 0 \rightarrow \mathcal{A}_n(k+1)/\mathcal{A}_n'(k+1) \rightarrow \mathcal{A}_n(k)/\mathcal{A}_n'(k+1)
     \rightarrow \mathrm{gr}^k (\mathcal{A}_n) \rightarrow 0 \]
induced from natural homomorphisms.

\subsection{Representation theory of $\mathrm{GL}(n,\Q)$}\label{Ss-Rep}
\hspace*{\fill}\ 

\vspace{0.5em}

Here we briefly review well-known results in representation theory for the general linear group $\mathrm{GL}(n,\Q)$,
including Cartan-Weyl's highest weight theory. The notation we use here is according to our previous paper \cite{ES1} by Enomoto-Satoh.

\vspace{0.5em}

First, we fix a basis $e_1, \ldots, e_n \in H_{\Q}$,
and by using it, we identify $\mathrm{GL}(H_{\Q})$ with $\mathrm{GL}(n,\Q)$.
Let $e_1^*, \ldots, e_n^*$ be its dual basis in $H_{\Q}^* := \mathrm{Hom}_{\Q}(H_{\Q}, \Q)$.
Let
\[ T_n := \{ \mathrm{diag}(t_1, \ldots ,t_n) \in \mathrm{GL}(n,\Q) \,|\, t_j \neq 0, \ 1 \leq j \leq n \} \]
be the maximal torus of $\mathrm{GL}(n,\Q)$.
For any $1 \leq i \leq n$,
define the one-dimensional representations $\varepsilon_i : T_n \rightarrow \Q^{\times}$ of $T_n$
by $\varepsilon_i(\mathrm{diag}(t_1, \ldots ,t_n))=t_i$.
Then 
\[\begin{split}
P    & := \{ \lambda_1\varepsilon_1+ \cdots +\lambda_n\varepsilon_n \ | \ \lambda_i \in \Z, \,\, 1 \leq i \leq n \}\cong \Z^n, \\
P^+  & := \{ \lambda_1\varepsilon_1+ \cdots +\lambda_n\varepsilon_n \in P \ | \ \lambda_1 \geq \lambda_2 \geq \cdots \geq \lambda_n\}
\end{split}\]
give the weight lattice and the set of dominant integral weights of $\mathrm{GL}(n,\Q)$ respectively. 
For simplicity, we write $\lambda=(\lambda_1, \ldots, \lambda_n) \in \Z^n$ for $\lambda= \lambda_1\varepsilon_1+ \cdots +\lambda_n\varepsilon_n \in P$ or $P^+$
if there is no confusion.

\vspace{0.5em}

For a rational representation $V$ of $\mathrm{GL}(n,\Q)$, consider the irreducible decomposition
$V=\bigoplus_{\lambda \in P}V_\lambda$ as a $T$-module where
\[ V_\lambda:=\{v \in V \ | \ tv=t_1^{\lambda_1} \cdots t_n^{\lambda_n} v \ \text{for any} \ t \in T \}.\]
We call this decomposition the weight decomposition of $V$ with respect to $T$.
If $V_\lambda \neq \{0\}$, then we call $\lambda$ the weight of $V$. For a weight $\lambda$,
a non-zero vector $v \in V_\lambda$ is call a weight vector of weight $\lambda$.

\vspace{0.5em}

Let $U$ be the subgroup of $\mathrm{GL}(n,\Q)$ consisting of all upper unitriangular matrices in $\mathrm{GL}(n,\Q)$.
For a rational representation $V$ of $\mathrm{GL}(n,\Q)$,
we set
\[ V^U := \{v \in V \ | \ uv=v \ \text{for all} \ u \in U\}. \]
We call a non-zero vector $v \in V^U$ a maximal vector of $V$.
Since $V^U$ is $T$-invariant subspace, we have the irreducible decomposition
$V^U=\bigoplus_{\lambda \in P}V^U_\lambda$ as a $T$-module where $V^U_\lambda:=V^U \cap V_\lambda$.
Then we have the following.
\begin{thm}[Cartan-Weyl's highest weight theory] \quad 
\begin{enumerate}
\item Any rational representation of $V$ is completely reducible.
\item Suppose $V$ is an irreducible rational representation of $\mathrm{GL}(n,\Q)$.
Then $V^U$ is one-dimensional, and the weight $\lambda$ of $V^U=V_\lambda^U$ belongs to $P^+$.
We call this $\lambda$ the highest weight of $V$, and any non-zero vector $v \in V^U_\lambda$ is called a highest weight vector of $V$.
\item For any $\lambda \in P^+$, there exists a unique (up to isomorphism) irreducible rational representation $L^\lambda$ of $\mathrm{GL}(n,\Q)$
with highest weight $\lambda$. Moreover, for two $\lambda, \mu \in P^+$, we have $L^\lambda \cong L^\mu$ if and only if $\lambda=\mu$.
\item The set of isomorphism classes of irreducible rational representations of $\mathrm{GL}(n,\Q)$ is parameterized by the set $P^+$ of dominant integral weights.
\item Let $V$ be a rational representation of $\mathrm{GL}(n,\Q)$ and $\chi_V$ the character of $V$ as a $T$-module.
Then for two rational representation $V$ and $W$, they are isomorphic as $\mathrm{GL}(n,\Q)$-modules if and only if $\chi_V=\chi_W$.
\end{enumerate}
\end{thm}

From the above theorem,
we can parameterize the set of isomorphism classes of irreducible rational representations of $\mathrm{GL}(n,\Q)$
by $P^+$. We can do this with the determinant representations. For any $e \in \Z$, let
$D^e:\mathrm{GL}(n,\Q) \rightarrow \Q^\times$ be the $e$-th power of the determinant representation of $\mathrm{GL}(n,\Q)$
defined by $X \mapsto (\det{X})^e$.
The highest weight of this representation is given by $(e,e, \cdots ,e) \in P^+$.
If $\lambda \in P^+$ satisfies $\lambda_n<0$ then we have
then
\[ L^\lambda \cong D^{-\lambda_n} \otimes_{\Q} L^{(\lambda_1-\lambda_n,\lambda_2-\lambda_n, \ldots, \lambda_{n-1} - \lambda_n,0)}. \]
Therefore we can parameterize the set of isomorphism classes of irreducible rational representations of $\mathrm{GL}(n,\Q)$
by the set $\{(\lambda,e)\}$ where $\lambda$ is a partition such that $\ell(\lambda) \leq n$ and $e \in \Z_{<0}$
where $\ell(\lambda)$ is the length of $\lambda$.
Moreover, the set of isomorphism classes of irreducible polynomial representations is parameterized by the set of partitions $\lambda$
such that $\ell(\lambda) \leq n$.
We remark that the dual representation of $L^{(\lambda_1,\lambda_2, \ldots ,\lambda_n)}$ is isomorphic to $L^{(-\lambda_n,\ldots ,-\lambda_2,-\lambda_1)}$.
In the following, for simplicity, for any $\lambda=(\lambda_1, \ldots, \lambda_n)$, we write the irreducible representation
$L^{\lambda}$ as $[\lambda_1, \ldots, \lambda_n]$ according to the usual notation in representation theory.

\vspace{0.5em}

Here we give a few examples.
The standard representation $H_{\Q}$ of $\mathrm{GL}(n,\Q)$ and its dual representation $H_{\Q}^*$
are irreducible representations with highest weight $[1,0, \ldots ,0]$ and $[0, \ldots ,0,-1]$ respectively. 
Hence, we have $H_{\Q} \cong [1]$ and $H_{\Q}^* \cong D^{-1} \otimes [1^{n-1}]$.
Set $W:= H_1(\mathrm{IA}_n,\Z) \cong H^* \otimes_{\Z} \Lambda^2 H$.
From Pieri's rule, the irreducible decomposition of $W_{\Q}$ as a $\mathrm{GL}(n,\Q)$-module is given by
\[ W_{\Q} \cong [1] \oplus (D^{-1} \otimes_{\Q} [2^2,1^{n-3}]). \]

\vspace{0.5em}

In Section {\rmfamily \ref{S-IA3}}, we consider the irreducible decompositions of several $\mathrm{GL}(W_{\Q})$-modules as a $\mathrm{GL}(H_{\Q})$-module
for $n=3$. Note that we identify $\mathrm{GL}(H_{\Q})$ with $\mathrm{GL}(n,\Q)$.
In order to find maximal vectors of them, we have to know the actions of elementary matrices of $\mathrm{GL}(n,\Q)$ on
the (coset classes of) Magnus generators in $W_{\Q}$. For any $1 \leq i \neq j \leq n$, let $E_{x_i x_j} \in \mathrm{Aut}\,F_n$ be the
Nielsen automorphism
of $F_n$ defined by
\[ x_p \mapsto \begin{cases}
                 x_i x_j, \hspace{1em} & p=i, \\
                 x_p, & p \neq i.
                \end{cases} \]
Let $E_{ij} \in \mathrm{GL}(n,\Z) \subset \mathrm{GL}(n,\Q)$ be the image of $E_{x_i x_j}$ by the natural homomorphism
$\rho : \mathrm{Aut}\,F_n \rightarrow \mathrm{GL}(n,\Z)$. Then the actions of $E_{ij}$ on the Magnus generators are given as follows:

\vspace{0.5em}

{\footnotesize
\begin{center}
{\renewcommand{\arraystretch}{1.4}
\begin{tabular}{|l|l|} \hline
  $K_{ij}^{E_{pq}} = K_{ij}$ & $K_{ijl}^{E_{pq}} = K_{ijl}$  \\ \hline
  $K_{ij}^{E_{il}} = K_{ij}+K_{ilj}$  & $K_{ijl}^{E_{iq}} = K_{ijl}$  \\ \hline
  $K_{ij}^{E_{jl}} = K_{ij}+ K_{il}$  & $K_{ijl}^{E_{jq}} = K_{ijl}+K_{iql}$  \\ \hline
  $K_{ij}^{E_{li}} = K_{ij}-K_{lij}$  & $K_{ijl}^{E_{pi}} = K_{ijl}+K_{plj}$  \\ \hline
  $K_{ij}^{E_{lj}} = K_{ij}$          & $K_{ijl}^{E_{pj}} = K_{ijl}$  \\ \hline
  $K_{ij}^{E_{ij}} = K_{ij}$          & $K_{ijl}^{E_{ij}} = K_{ijl}$  \\ \hline
  $K_{ij}^{E_{ji}} = K_{ij}+K_{ji}$   & $K_{ijl}^{E_{il}} = K_{ijl}$  \\ \hline
  & $K_{ijl}^{E_{ji}} = K_{ijl}+K_{jli}+K_{il}-K_{jl}$ \\ \hline
  & $K_{ijl}^{E_{jl}} = K_{ijl}$ \\ \hline
  & $K_{ijl}^{E_{lj}} = K_{ijl}$ \\ \hline
\end{tabular}}
\end{center}
}

\vspace{0.5em}

\noindent
Here $1 \leq i, j, l, p, q \leq n$ are distinct indices.

\section{(Co)homology groups of $\mathrm{IA}_n$}\label{S-sec}

\vspace{0.5em}

In this section, we give combinatorial group theoretic descriptions of low dimensional (co)homology groups of $\mathrm{IA}_n$.
In particular, by using it, we study the second and the third (co)homology groups of $\mathrm{IA}_3$.

\vspace{0.5em}

Since $\mathrm{IA}_3$ is not finitely presentable due to Krsti\'{c} and McCool \cite{Krs}, there is a possibility that $H_2(\mathrm{IA}_3, \Z)$
is not finitely generated. Bestvina, Bux and D. Margalit \cite{Bes} showed that $H_2(\mathrm{IA}_3, \Z)$ is not finitely generated.
It is a natural problem to give an explicit generating set of $H_2(\mathrm{IA}_3,\Z)$.
However, if we approach this problem, we immediately face the difficulties coming from the complexity of the structure of the group
of relators among Magnus generators of $\mathrm{IA}_3$.
In order to make the problem more easy, we consider the rationalization of the problem. After this, we can use the representation theory of the general linear group
$\mathrm{GL}(n,\Q)$, and its highest weight theory.
But, it is still too complicated to give a complete answer to the above problem. In the following, we give partial results for this problem.

\vspace{0.5em}

\subsection{A minimal presentation and $H_2(\mathrm{IA}_n,\Z)$}\label{Ss-min}
\hspace*{\fill}\ 

\vspace{0.5em}

In this subsection, we give a combinatorial group theoretical interpretation of the second homology group of $\mathrm{IA}_n$.
Let $F$ be a free group on the Magnus generators: $K_{ij}$ for any $1 \leq i \neq j \leq n$, and $K_{ijl}$ for any $i \neq j<l \neq i$.
By abuse of the language, we use the same notation $K_{ij}$ and $K_{ijl}$ for the elements in $F$
as the automorphisms $K_{ij}$ and $K_{ijl}$ in $\mathrm{IA}_n$.
The rank of $F$ is $n^2(n-1)/2$.
We have a natural surjective homomorphism $\pi : F \rightarrow \mathrm{IA}_n$, and the group extension
\begin{equation}\label{eq-ext}
1 \rightarrow R \rightarrow F \xrightarrow{\pi} \mathrm{IA}_n \rightarrow 1
\end{equation}
of $\mathrm{IA}_n$ where $R=\mathrm{Ker}(\pi)$. Since the abelianization
$H_1(\mathrm{IA}_n,\Z)$ of $\mathrm{IA}_n$ is the free abelian group generated by
(the coset classes of) the Magnus generators, we verify that $\pi$ induces the isomorphism
\[ \pi_* : H_1(F,\Z) \rightarrow H_1(\mathrm{IA}_n,\Z). \]
This shows that $R$ is contained in the commutator subgroup $[F,F]$ of $F$.
From the homological five-term exact sequence of (\ref{eq-ext}),
we have
\[\begin{split}
   H_2(F,\Z) \rightarrow H_2(\mathrm{IA}_n,\Z) \rightarrow H_1(R,\Z)_F 
    \rightarrow H_1(F,\Z) \xrightarrow{\pi_*} H_1(\mathrm{IA}_n,\Z) \rightarrow 0,
  \end{split}\]
and hence
\[ H_2(\mathrm{IA}_n,\Z) \cong H_1(R,\Z)_F. \]

\vspace{0.5em}

Let $F=\Gamma_F(1) \supset \Gamma_F(2) \supset \cdots$ be the lower central series of $F$, and set
$\mathcal{L}_F(k) := \Gamma_F(k)/\Gamma_F(k+1)$ for each $k \geq 1$.
Let $R = R_1 \supset R_2 \supset \cdots$ be the descending filtration of $R$ defined by $R_k := R \cap \Gamma_F(k)$ for each $k \geq 1$. 
We have $R_k=R$ for $k=1$ and $2$.
For each $k \geq 1$, let
\[ \pi_k : \mathcal{L}_F(k) \rightarrow \mathrm{gr}^k({\mathcal{A}_n'}) \]
be the homomorphism induced from the natural projection $\pi : F \rightarrow \mathrm{IA}_n$.
By observing $R_k/R_{k+1} \cong (R_k \, \Gamma_F(k+1))/\Gamma_F(k+1)$, we obtain the exact sequence
\begin{equation}\label{eq-momo}
0 \rightarrow R_k/R_{k+1} \rightarrow \mathcal{L}_F(k) \xrightarrow{\pi_k} \mathrm{gr}^k({\mathcal{A}_n'}) \rightarrow 0.
\end{equation}
For each $k \geq 2$.
The natural projection $R \rightarrow R/R_{k+1}$ induces the surjective homomorphism
\[ \psi_k : H_1(R,\Z) \rightarrow H_1(R/R_{k+1},\Z). \]
By considering the right action of $F$ on $R$, defined by
\[ r \cdot x := x^{-1} r x, \hspace{1em} r \in R, \,\,\, x \in F, \]
we see $\psi_k$ is an $F$-equivariant homomorphism.
Hence it induces the surjective homomorphism
\[ H_1(R,\Z)_{F} \rightarrow H_1(R/R_{k+1},\Z)_{F}, \]
which is also denoted by $\psi_k$.
For $k=2$, $H_1(R/R_{3},\Z)_{F} = R/R_{3}$ since $F$ acts on $R/R_{3}$ trivially.
Therefore we can detect non-trivial elements
in $H_2(\mathrm{IA}_n,\Z)$ by using $R/R_3$ if $R/R_{3}$ is a non-trivial.

\vspace{0.5em}

In the exact sequence (\ref{eq-momo}) for $k=2$,
since $\mathcal{L}_F(2)$ is a free abelian group,
so is $R/R_{3}$. Since $\mathcal{A}_n'(3) = \mathcal{A}_n(3)$ by \cite{S25},
we see $\mathrm{gr}^2(\mathcal{A}_n') = \mathrm{gr}^2(\mathcal{A}_n)$.
Hence, we have
\[\begin{split}
   \mathrm{rank}_{\Z}(R/R_{3})
    & = \mathrm{rank}_{\Z}(\mathcal{L}_F(2)) - \mathrm{rank}_{\Z}(\mathrm{gr}^2(\mathcal{A}_n)) \\
    & = \frac{1}{8}n^2(n-1)(n^3-n^2-2) - \frac{1}{6}n(2n^3-5n-3).
   \end{split}\]
We see that $H_2(\mathrm{IA}_n,\Z)$ contains a free abelian group of the above rank.

\vspace{0.5em}

Next, we consider the kernel of $\psi_2 : H_1(R,\Z)_{F} \rightarrow H_1(R/R_3,\Z)_{F}$.
Observe the exact sequences
\begin{eqnarray}
 0 \rightarrow R_3/[F,R] \rightarrow & R/[F,R] & \xrightarrow{\psi_2} R/R_3 \rightarrow 0, \label{eq-aqua} \\ 
 0 \rightarrow (R_4 [F,R])/[F,R] \rightarrow & R_3/[F,R] & \rightarrow R_3/(R_4 [F,R]) \rightarrow 0 \label{eq-mint}
\end{eqnarray}
of $\Z$-modules. This shows that if $H_2(\mathrm{IA}_n,\Z) \cong R/[F,R]$ is not isomorphic to $R/R_3$, by showing $R_3/(R_4 [F,R]) \neq 0$,
we have a potential to detect non-trivial second homology classes which we cannot be detected by $R/R_3$.
In Section {\rmfamily \ref{S-IA3}}, for $n=3$ we show that $(R_3/(R_4 [F,R])) \otimes_{\Z} \Q$ is a non-trivial irreducible representation of $\mathrm{GL}(3,\Q)$.
Consider the homomorphism
\[ \mathcal{L}_F(2) \otimes_{\Z} \mathcal{L}_F(1) \rightarrow \mathcal{L}_F(3) \]
induced from the commutator bracket of $F$. This homomorphism induces the homomorphism
\[ [\, , \,] : R/R_3 \otimes_{\Z} \mathcal{L}_F(1) \rightarrow R_3/R_4. \]
Then $R_3/(R_4 [F,R]) \cong \mathrm{Coker}([\, , \,])$.

\vspace{0.5em}

\subsection{On the second cohomology group $H^2(\mathrm{IA}_n,\Z)$}\label{Ss-coh}
\hspace*{\fill}\ 

\vspace{0.5em}

Here we consider the second cohomology group of $\mathrm{IA}_n$, and
the image of the cup product $\cup : \Lambda^2 H^1(\mathrm{IA}_n,\Z) \rightarrow H^2(\mathrm{IA}_n,\Z)$.
From the cohomological five-term exact sequence of (\ref{eq-ext}), we have
\[\begin{split}
    0 \rightarrow H^1(\mathrm{IA}_n,\Z) \xrightarrow{\pi^*} H^1(F,\Z)
    \rightarrow H^1(R,\Z)^F \rightarrow H^2(\mathrm{IA}_n,\Z) \rightarrow H^2(F,\Z),
  \end{split}\]
and hence
\[ H^2(\mathrm{IA}_n,\Z) \cong H^1(R,\Z)^F. \]
For any $k \geq 2$, the natural projection $R \rightarrow R/R_{k+1}$ induces the injective homomorphism
\[ \psi^k : H^1(R/R_{k+1},\Z)^F \rightarrow H^1(R,\Z)^F. \]
In particular, for $k=2$, we have $H^1(R/R_{3},\Z)^F \cong H^1(R/R_{3},\Z)$ since $F$ acts on $R/R_3$ trivially.
Namely we can regard $H^1(R/R_{3},\Z)$ as a $\Z$-submodule of $H^2(\mathrm{IA}_n,\Z)$.
Then we have
\begin{lem}\label{PP}
The image of the cup product
\[ \cup : \Lambda^2 H^1(\mathrm{IA}_n,\Z) \rightarrow H^2(\mathrm{IA}_n,\Z) \]
is isomorphic to the image of
\[ \iota^* : H^1(\mathcal{L}_F(2),\Z) \rightarrow H^1(R/R_3,\Z) \]
induced from the natural inclusion $R/R_3 \rightarrow \mathcal{L}_F(2)$.
\end{lem}
\textit{Proof.}
First, observing the last three terms of the cohomological five-term exact sequence of
\begin{equation}\label{tau1}
1 \rightarrow \mathcal{A}_n'(2) \rightarrow \mathrm{IA}_n \rightarrow \mathrm{IA}_n^{\mathrm{ab}} \rightarrow 1,
\end{equation}
we have
\[\begin{split}
     0 \rightarrow H^1(\mathcal{A}_n'(2),\Z)^{\mathrm{IA}_n}
          \rightarrow H^2(\mathrm{IA}_n^{\mathrm{ab}},\Z) \rightarrow H^2(\mathrm{IA}_n,\Z).
  \end{split}\]
Since $H^1(\mathcal{A}_n'(2),\Z)^{\mathrm{IA}_n}= H^1(\mathrm{gr}^2 (\mathcal{A}_n'),\Z)$,
we obtain an exact sequence
\[ 0 \rightarrow H^1(\mathrm{gr}^2 (\mathcal{A}_n'),\Z)
      \rightarrow H^2(\mathrm{IA}_n^{\mathrm{ab}},\Z) \xrightarrow{} H^2(\mathrm{IA}_n,\Z). \]
Since $\mathrm{IA}_n^{\mathrm{ab}}$ is a free abelian group of finite rank,
we have the natural isomorphism $H^2(\mathrm{IA}_n^{\mathrm{ab}},\Z) \cong \Lambda^2 H^1(\mathrm{IA}_n,\Z)$.
Then the map $H^2(\mathrm{IA}_n^{\mathrm{ab}},\Z) \xrightarrow{} H^2(\mathrm{IA}_n,\Z)$ is regarded as
the cup product $\cup : \Lambda^2 H^1(\mathrm{IA}_n, \Z) \rightarrow H^2(\mathrm{IA}_n,\Z)$.

\vspace{0.5em}

On the other hand, from the cohomological five-term exact sequence
\[\begin{split}
     0 \rightarrow H^1(\mathrm{gr}^2 & (\mathcal{A}_n),\Z) \rightarrow H^1(\mathcal{L}_F(2), \Z)
       \rightarrow H^1(R/R_3,\Z)^{\mathcal{L}_F(2)} \\
         & \rightarrow H^2(\mathrm{gr}^2(\mathcal{A}_n'),\Z) \rightarrow H^2(\mathcal{L}_F(2),\Z)
  \end{split}\]
of (\ref{eq-momo}) for $k=2$,
since $\mathcal{L}_F(2)$ acts on $R/R_3$ trivially, we have $H^1(R/R_3,\Z)^{\mathcal{L}_F(2)}
 = H^1(R/R_3,\Z)$.
Then we have the commutative diagram
\[\begin{CD}
     0     @>>>  H^1(\mathrm{gr}^2 (\mathcal{A}_n'),\Z) @>{\mathrm{tg}}>> H^2(\mathrm{IA}_n^{\mathrm{ab}},\Z)  @>{\cup}>>  H^2(\mathrm{IA}_n,\Z)     \\
     @.               @|                                    @VV{\mu}V                    @.          \\
     0     @>>>  H^1(\mathrm{gr}^2 (\mathcal{A}_n'),\Z) @>>> H^1(\mathcal{L}_F(2),\Z) @>{\iota^*}>>
     H^1(R/R_3,\Z)
 \end{CD}\]
where $\mathrm{tg}$ is the transgression, $\mu$ is the natural isomorphism, and $\iota^*$ is the homomorphism
induced from the inclusion $R/R_3 \rightarrow \mathcal{L}_F(2)$.
Hence we obtain $\mathrm{Im}(\cup) \cong \mathrm{Im}(\iota^*)$. $\square$

\vspace{0.5em}

We remark that since $\mathcal{A}_n'(3) = \mathcal{A}_n(3)$, and since $\mathrm{gr}^2({\mathcal{A}_n'}) = \mathrm{gr}^2({\mathcal{A}_n})$ is a free abelian group,
the map $\iota^*$ is surjective. This implies that $\mathrm{Im}(\cup) \cong H^1(R/{R}_3,\Z)$.

\vspace{0.5em}

Next, we consider a method to detect non-trivial elements in $H^2(\mathrm{IA}_n,\Z)$, which do not belong to $\mathrm{Im}(\cup)$.
This is a cohomology version of the argument in Subsection {\rmfamily \ref{Ss-min}}.
From the exact sequence
\[ 0 \rightarrow R_4[R,F]/[R,F] \rightarrow R/[F,R] \rightarrow R/R_4[R,F] \rightarrow 0, \]
we have the exact sequence
\[ 0 \rightarrow \mathrm{Hom}_{\Z}(R/R_4[R,F],\Z) \rightarrow \mathrm{Hom}_{\Z}(R/[F,R],\Z) \rightarrow \mathrm{Hom}_{\Z}(R_4[R,F]/[R,F],\Z). \]
Here $\mathrm{Hom}_{\Z}(R/[F,R],\Z) \cong H^2(\mathrm{IA}_n,\Z)$.
On the other hand, since $R/R_3$ is a free abelian group, the split exact sequence
\[ 0 \rightarrow R_3/R_4[R,F] \rightarrow R/R_4[F,R] \rightarrow R/R_3 \rightarrow 0 \]
as a $\Z$-module induces the exact sequence
\[ 0 \rightarrow \mathrm{Hom}_{\Z}(R/R_3,\Z) \rightarrow \mathrm{Hom}_{\Z}(R/R_4[F,R],\Z) \rightarrow \mathrm{Hom}_{\Z}(R_3/R_4[R,F],\Z)
     \rightarrow 0. \]
Thus by studying $\mathrm{Hom}_{\Z}(R_3/R_4[R,F],\Z)$, we have a potential to detect a non-trivial submodule in
$H^2(\mathrm{IA}_n,\Z)$ which is not contained in $\mathrm{Im}(\cup)$.

\vspace{0.5em}

\subsection{On the third cohomology group $H^3(\mathrm{IA}_n,\Q)$}\label{Ss-coh}
\hspace*{\fill}\ 

\vspace{0.5em}

Here we consider the rational third cohomology group of $\mathrm{IA}_n$. In particular, we characterize the image of the triple cup product
\[ \cup_{\mathrm{IA}_n}^3 : \Lambda^3 H^1(\mathrm{IA}_n, \Q) \rightarrow H^3(\mathrm{IA}_n, \Q). \]

\vspace{0.5em}

In general, by the reduction theorem, we have the isomorphism
\[ H^1(\mathrm{IA}_n, \mathrm{Hom}_{\Q}(R^{\mathrm{ab}},\Q)) \xrightarrow{\cong} H^3(\mathrm{IA}_n, \Q). \]
This isomorphism is given by the differential homomorphism $d_2^{1,1}$ of the Lyndon-Hochshild-Serre spectral sequence of
the group extension (\ref{eq-ext}). (See \cite{Hil} for details.)
The $\mathrm{IA}_n$-equivariant exact sequence
\[ 0 \rightarrow R_3/[R,R] \xrightarrow{\alpha} R^{\mathrm{ab}} \xrightarrow{\beta} R/R_3 \rightarrow 0 \]
induces the $\mathrm{IA}_n$-equivariant exact sequence
\[ 0 \rightarrow \mathrm{Hom}_{\Q}(R/R_3,\Q) \xrightarrow{\beta^*} \mathrm{Hom}_{\Q}(R^{\mathrm{ab}},\Q)
     \xrightarrow{\alpha^*} \mathrm{Hom}_{\Q}(R_3/[R,R],\Q) \rightarrow 0. \]
In the following, for simplicity, for any $\Q$-vector space $V$, we write $V_{\Q}^*$ for the $\Q$-linear dual space $\mathrm{Hom}_{\Q}(V,\Q)$
of $V$. Then the above exact sequence induces the cohomological long exact sequence
\begin{eqnarray*}
 0 & \rightarrow & H^0(\mathrm{IA}_n, (R/R_3)_{\Q}^*) \xrightarrow{\beta^0} H^0(\mathrm{IA}_n, (R^{\mathrm{ab}})_{\Q}^*)
       \xrightarrow{\alpha^0} H^0(\mathrm{IA}_n, (R_3/[R,R])_{\Q}^*) \\ 
   & \xrightarrow{\delta} & H^1(\mathrm{IA}_n, (R/R_3)_{\Q}^*) \xrightarrow{\beta^1} H^1(\mathrm{IA}_n, (R^{\mathrm{ab}})_{\Q}^*)
       \xrightarrow{\alpha^1}  H^1(\mathrm{IA}_n, (R_3/[R,R])_{\Q}^*) \\
   & \rightarrow & \cdots
\end{eqnarray*}
of $\mathrm{IA}_n$. Since $\mathrm{IA}_n$ acts on $R_3/R_4$ and $R_3/R_4[F,R]$ trivially, we have
the sequence
\[ (R_3/R_4[F,R])_{\Q}^* \hookrightarrow (R_3/R_4)_{\Q}^* \hookrightarrow H^0(\mathrm{IA}_n, (R_3/[R,R])_{\Q}^*) \]
of natural $\Q$-linear injective homomorphisms. Through these maps, we consider $(R_3/R_4[F,R])_{\Q}^*$ and $(R_3/R_4)_{\Q}^*$ as
submodules of $H^0(\mathrm{IA}_n, (R_3/[R,R])_{\Q}^*)$.

\vspace{0.5em}

Now we consider the restriction of the connecting homomorphism $\delta$ to $(R_3/R_4)_{\Q}^*$. To begin with, we determine
the kernel of $\delta$.
\begin{lem}\label{L-conn}
$\mathrm{Im}(\alpha^0) \cap (R_3/R_4)_{\Q}^* = (R_3/R_4[F,R])_{\Q}^*$.
\end{lem}
\textit{Proof}.
Consider the natural isomorphisms
\begin{eqnarray*}
 H^0(\mathrm{IA}_n, (R^{\mathrm{ab}})_{\Q}^*) & \cong & \mathrm{Hom}_{\Q}(R \big{/} [F,R], \Q), \\
 H^0(\mathrm{IA}_n, (R_3/[R,R])_{\Q}^*) & \cong & \mathrm{Hom}_{\Q}(R_3 \big{/} [R,R][R_3,F], \Q). 
\end{eqnarray*}
Then $\alpha^0 : \mathrm{Hom}_{\Q}(R/[F,R], \Q) \rightarrow \mathrm{Hom}_{\Q}(R_3/[R,R][R_3,F], \Q)$ is given by
\[ f \mapsto (r \pmod{[R,R][R_3,F]} \mapsto f(r \pmod{[F,R]})). \]
On the other hand, the exact sequence
\[ 0 \rightarrow R_3/R_4[R,R] \rightarrow R/R_4[F,R] \xrightarrow R/R_3 \rightarrow 0 \]
splits since $R/R_3$ is a free abelian group of finite rank, we have an isomorphism
\[ R/R_4[F,R] \cong R_3/R_4[R,R] \oplus R/R_3 \]
as $\Z$-modules, and hence $(R/R_4[F,R])_{\Q} \cong (R_3/R_4[R,R])_{\Q} \oplus (R/R_3)_{\Q}$ as $\Q$-vector spaces.
Fix this isomorphism, and let $p : (R/R_4[F,R])_{\Q} \rightarrow (R_3/R_4[R,R])_{\Q}$ be the projection map.

\vspace{0.5em}

For any $g \in (R_3/R_4[F,R])_{\Q}^*$, consider the extension $\tilde{g} := g \circ p : (R/R_4[F,R])_{\Q} \rightarrow \Q$
of $g$. Set $\overline{g} := \tilde{g} \circ \gamma_{\Q} : (R/[F,R])_{\Q} \rightarrow \Q$ where $\gamma : R/[F,R] \rightarrow R/R_4[F,R]$
is the natural map. Then we see $\alpha^0(\overline{g})=g$.

\vspace{0.5em}

Conversely, for any $h \in \mathrm{Im}(\alpha^0) \cap (R_3/R_4)_{\Q}^*$, $h$ vanishes on $R_4[F,R]/R_4$. Thus we see
$h \in (R_3/R_4[F,R])_{\Q}^* \subset H^0(\mathrm{IA}_n, (R_3/[R,R])_{\Q}^*)$.
$\square$

\vspace{0.5em}

Next we consider the image of $\alpha^0$ in $H^1(\mathrm{IA}_n, (R^{\mathrm{ab}})_{\Q}^*)$.
Observe the commutative diagram
\[\begin{CD}
     0     @>>>  R_3/[R,R] @>>> R/[R,R]  @>>>  R/R_3   @>>> 0  \\
     @.               @VVV                                    @VVV          @VVV          @.          \\
     0     @>>>  \Gamma_F(3)/[\Gamma_F(2), \Gamma_F(2)] @>>> \Gamma_F(2)^{\mathrm{ab}} @>>> \Gamma_F(2)/\Gamma_F(3) @>>> 0
 \end{CD}\]
of $\mathrm{IA}_n^{\mathrm{ab}}$-modules, containing two exact sequences.
This diagram and the natural projection $\mathrm{IA}_n \rightarrow \mathrm{IA}_n^{\mathrm{ab}}$ induces the commutative diagram
\begin{equation}\label{eq-peace}\begin{CD}
     H^1(\mathrm{IA}_n^{\mathrm{ab}}, (\Gamma_F(2)/\Gamma_F(3))_{\Q}^*)   @>{\widetilde{\beta}^1}>> 
          H^1(\mathrm{IA}_n^{\mathrm{ab}}, (\Gamma_F(2)^{\mathrm{ab}})_{\Q}^*) @>{d_2^{1,1}}>> H_3(\mathrm{IA}_n^{\mathrm{ab}},\Q) \\
               @V{\psi_1}VV                                    @VV{\xi_1}V         \\
     H^1(\mathrm{IA}_n^{\mathrm{ab}}, (R/R_3)_{\Q}^*)   @>>> 
          H^1(\mathrm{IA}_n^{\mathrm{ab}}, (R^{\mathrm{ab}})_{\Q}^*)  \\
               @V{\psi_2}VV                                    @VV{\xi_2}V         \\
     H^1(\mathrm{IA}_n, (R/R_3)_{\Q}^*)   @>{\beta^1}>> 
          H^1(\mathrm{IA}_n, (R^{\mathrm{ab}})_{\Q}^*) @>{d_2^{1,1}}>> H_3(\mathrm{IA}_n,\Q)  \\
 \end{CD}\end{equation}
of the first cohomologies of $\mathrm{IA}_n$ and $\mathrm{IA}_n^{\mathrm{ab}}$. Here $\psi_1$ is surjective and $\psi_2$ is an isomorphism.

\vspace{0.5em}

Remark that through the isomorphism $d_2^{1,1}$, the target of $\widetilde{\beta}^1$ is isomorphic to the image of the triple cup product
\[ \cup_{\mathrm{IA}_n^{\mathrm{ab}}}^3 : \Lambda^3 H^1(\mathrm{IA}_n^{\mathrm{ab}}, \Q) \rightarrow H^3(\mathrm{IA}_n^{\mathrm{ab}}, \Q) \]
since $\cup_{\mathrm{IA}_n^{\mathrm{ab}}}^3$ is isomorphism.

\vspace{0.5em}

Now we show the main proposition of this subsection.

\begin{pro}\label{P-cup3-2}
Through the isomorphism $d_2^{1,1}$, the image of $\beta^1$ is isomorphic to the image of the triple cup product
\[ \cup_{\mathrm{IA}_n}^3 : \Lambda^3 H^1(\mathrm{IA}_n, \Q) \rightarrow H^3(\mathrm{IA}_n, \Q). \]
\end{pro}
\textit{Proof}.
For any $f \in \mathrm{Im}(\beta^1)$, from the surjectivity of $\psi_2 \circ \psi_1$, there exists some element
$\widetilde{f} \in H^1(\mathrm{IA}_n^{\mathrm{ab}}, (\Gamma_F(2)/\Gamma_F(3))_{\Q}^*)$ such that
$(\beta^1 \circ (\psi_2 \circ \psi_1))(\widetilde{f})=f$. By the commutativity of (\ref{eq-peace}), we see
$\widetilde{\beta}^1 \circ (\xi_2 \circ \xi_1)(\widetilde{f})=f$.
With standard arguments and calculations, we can see that the homomorphism
$\xi_2 \circ \xi_1 : H^3(\mathrm{IA}_n^{\mathrm{ab}}, \Q) \rightarrow H^3(\mathrm{IA}_n, \Q)$
is equal to the natural homomorphism induced from that projection $\mathrm{IA}_n \rightarrow \mathrm{IA}_n^{\mathrm{ab}}$.
Therefore, we see $f \in \mathrm{Im}(\cup_{\mathrm{IA}_n}^3)$.

\vspace{0.5em}

On the other hand, since $H^1(\mathrm{IA}_n,\Q) \cong H^1(\mathrm{IA}_n^{\mathrm{ab}},\Q)$, for any
$g \in \mathrm{Im}(\cup_{\mathrm{IA}_n}^3)$, there exists some $\widetilde{g} \in \mathrm{Im}(\cup_{\mathrm{IA}_n^{\mathrm{ab}}}^3)$
such that $\xi_2 \circ \xi_1(\widetilde{g}) =g$. Hence there exists some $h \in H^1(\mathrm{IA}_n^{\mathrm{ab}}, (\Gamma_F(2)/\Gamma_F(3))_{\Q}^*)$
such that $\widetilde{\beta}^1 (h)=\widetilde{g}$ and
\[ {\beta}^1 \circ (\psi_2 \circ \psi_1) (h) = (\xi_2 \circ \xi_1) \circ \widetilde{\beta}^1 (h) = g. \]
Thus $g \in \mathrm{Im}(\beta^1)$, and we obtain the required result. $\square$

\vspace{0.5em}

\section{The case where $n=3$}\label{S-IA3}
\hspace*{\fill}\ 

\vspace{0.5em}

In this section, we consider the case where $n=3$ with rational representations of $\mathrm{GL}(3,\Q)$,
and show our main theorems.

\vspace{0.5em}

To begin with, we remark that for the standard $\mathrm{GL}(n,\Z)$-representation $H$, we regard $H_{\Q}$ as a $\mathrm{GL}(n,\Q)$-module in a usual way.
For the $\mathrm{GL}(n,\Z)$-module $W=H_1(\mathrm{IA}_n,\Z)$, we regard $W_{\Q}=H_1(\mathrm{IA}_n,\Q)=H_{\Q}^* \otimes_{\Q} \Lambda^2 H_{\Q}$ as a
$\mathrm{GL}(n,\Q)$-module via the standard representation $H_{\Q}$.
Since the rational free Lie algebra $:= \bigoplus_{k \geq 1} \mathcal{L}_F^{\Q}(k)$ is the free Lie algebra generated by $W_{\Q}$,
we can consider each of $\mathcal{L}_F^{\Q}(k)$
as a $\mathrm{GL}(W_{\Q})$-module, and hence $\mathrm{GL}(n,\Q)$-module.
In this section, we consider the irreducible decompositions of several submodules of $\mathcal{L}_F^{\Q}(k)$ as a $\mathrm{GL}(n,\Q)$-module, and
give maximal vectors corresponding to each of irreducible component. In the following, we always assume $n=3$.

\vspace{0.5em}

First, the irreducible decomposition of $W_{\Q}$ is given by
\[ W_{\Q} = [1] \oplus [1,1,-1]. \]
The list of maximal vectors of $W_{\Q}$ is given as follows:

\vspace{0.5em}

\begin{center}
{\renewcommand{\arraystretch}{1.2}
\begin{tabular}{|l|l|l|} \hline
  $V$  & maximal vectors & $\mathrm{dim}_{\Q} V$ \\ \hline
  $[1]$        & $\iota_1$     & $3$   \\ \hline
  $[1,1,-1] $  & $K_{312}$     & $6$   \\ \hline
\end{tabular}}
\end{center}

\vspace{0.5em}

Here we recall the definition of $\iota_j$. 
For any $1 \leq j \leq n$, $\iota_j$ is the inner automorphism of $F_n$ given by $x \mapsto x_j^{-1} x x_j$.
Hence, its coset class in $W_{\Q}$, also denoted by $\iota_j$, is given by 
$\iota_j=K_{1j}+K_{2j} +\cdots + K_{3j}$ for any $1 \leq j \leq n$.

\vspace{0.5em}

In \cite{Pet}, Pettet gave the irreducible decomposition of $\mathcal{L}_F^{\Q}(2)=\Lambda^2 W_{\Q}$ by
\begin{equation}\label{eq-pet}
 \Lambda^2 W_{\Q} \cong [1,1]^{\oplus 2} \oplus [2,1,-1]^{\oplus 2}.
\end{equation}
The list of linearly independent maximal vectors corresponding to this decomposition are as follows:

\vspace{0.5em}

\begin{center}
{\renewcommand{\arraystretch}{1.2}
\begin{tabular}{|l|l|l|} \hline
  $V$            & maximal vectors           & $\mathrm{dim}_{\Q} V$   \\ \hline
  $[1,1]$        & $\bm{v}_1:=[\iota_1, \iota_2]$  & $3$   \\
                 & $\bm{v}_2:=[K_{312},\iota_3]+[\iota_1, K_{12}] -2[K_{32},\iota_1] +[K_{31}, \iota_2]$ & \\ \hline
  $[2,1,-1]$     & $\bm{v}_3:=[K_{312}, K_{21}]$  & $15$  \\
                 & $\bm{v}_4:=[K_{312},\iota_1]$  &   \\ \hline
\end{tabular}}
\end{center}

\vspace{0.5em}

For any $k \geq 1$, consider the injective homomorphism $(R_k/R_{k+1})_{\Q} \rightarrow \mathcal{L}^{\Q}_F(k)$ induced from the natural inclusion map
$R_k \hookrightarrow \Gamma_F(k)$. We regard $(R_k/R_{k+1})_{\Q}$ as a $\mathrm{GL}(3,\Q)$-submodule of $\mathcal{L}^{\Q}_F(k)$.

\begin{lem}\label{L-mare}
$(R/R_3)_{\Q} \cong [1,1] \oplus [2,1,-1]$.
\end{lem}
\textit{Proof}.
First, we check $\bm{v}_2, \bm{v}_4 \in (R/R_3)_{\Q}$.
Observe the exact sequence (\ref{eq-momo}) for $k=2$, and recall that the second Johnson homomorphism
$\tau_{2,\Q}' : \mathrm{gr}^2_{\Q}(\mathcal{A}_3') \rightarrow H^*_{\Q} \otimes_{\Q} \mathcal{L}_3^{\Q}(3)$ is injective
since $\mathcal{A}_3'(k)=\mathcal{A}_3(k)$ for $k=2$ and $3$.
Thus, from $\tau_{2,\Q}' \circ \pi_{2,\Q}(\bm{v}_i)=0$ for $i=2$ and $4$, we see $\bm{v}_2, \bm{v}_4 \in (R/R_3)_{\Q}$.

\vspace{0.5em}

This shows that $(R/R_3)_{\Q}$ contains a submodule isomorphic to $[1,1] \oplus [2,1,-1]$.
On the other hand, from the fact that
\[ \mathrm{dim}_{\Q}((R/R_3)_{\Q}) = \mathrm{dim}_{\Q}(\mathcal{L}_F^{\Q}(2)) - \mathrm{dim}_{\Q}(\mathrm{gr}^2_{\Q}(\mathcal{A}_3'))
   = 36-18=18 \]
and $\mathrm{dim}_{\Q}([1,1] \oplus [2,1,-1])=18$, we obtain the required results. $\square$

\vspace{0.5em}

Next we give the irreducible decomposition of $\mathcal{L}_F^{\Q}(3)$. To do this, we prepare some lemmas.
In general, for any rational representations $V_1$, $V_2$ of $\mathrm{GL}(n,\Q)$, and for any $k \geq 1$ we have the formula
\begin{equation}\label{eq-ext-dec}
\Lambda^k (V_1 \oplus V_2) \cong \bigoplus_{p+q=k} (\Lambda^p V_1 \otimes_{\Q} \Lambda^q V_2).
\end{equation}
\noindent
Then we have the following decompositions.

\begin{lem}\label{L-aoi} \,\ \\
\vspace{-1em}

\begin{enumerate}
\item $\Lambda^2 [1,1,-1] \cong [2,1,-1]$.
\item $[3,2] \otimes_{\Q} [2,2] \cong [5,2,2] \oplus [4,3,2] \oplus [5,3,1] \oplus [5,4] \oplus [4,4,1]$.
\item $\Lambda^3 [1,1,-1] \cong [3] \oplus [2,2,-1]$.
\end{enumerate}
\end{lem}
\textit{Proof}.
(1). \,\,By applying (\ref{eq-ext-dec}) to $\Lambda^2 W_{\Q} = \Lambda^2 ([1] \oplus [1,1,-1])$, we have
\[ \Lambda^2 W_{\Q} \cong \Lambda^2 [1] \bigoplus ( [1] \otimes_{\Q} [1,1,-1] ) \bigoplus \Lambda^2 [1,1,-1]. \]
Since $\Lambda^2 [1] = [1,1]$ and since $[1] \otimes_{\Q} [1,1,-1] \cong [1,1] \oplus [2,1,-1]$ by Pieri's rule, we obtain
the required decomposition by (\ref{eq-pet}).

\vspace{0.5em}

(2). \,\,Since $n=3$, we have
\[ \Lambda^2 [1] \otimes_{\Q} \Lambda^2 [1] \cong (D \otimes_{\Q} [1] ) \bigoplus [2,2] \]
by Pieri's rule. On the other hand, by using Pieri's rule again, we can compute the irreducible decomposition of both of
$(\Lambda^2 [1] \otimes_{\Q} \Lambda^2 [1]) \otimes_{\Q} [3,2]$ and $(D \otimes_{\Q} [1] ) \otimes_{\Q} [3,2]$.
Hence by checking the difference between them, we obtain the required results.

\vspace{0.5em}

(3). \,\,The module $\Lambda^3 [1,1,-1]$ is a $\mathrm{GL}(n,\Q)$-submodule of
\[\begin{split}
   \Lambda^2 [1,1,-1] & \otimes_{\Q} [1,1,-1] \cong D^{-2} \otimes_{\Q} ([3,2] \otimes_{\Q} [2,2]) \\
   & = [3] \oplus [2,1] \oplus [3,1,-1] \oplus [3,2,-2] \oplus [2,2,-1].
  \end{split}\]
Since $\mathrm{dim}_{\Q} \Lambda^3 [1,1,-1]=20$, the unique possibility is $\Lambda^3 [1,1,-1] \cong [3] \oplus [2,2,-1]$.
$\square$

\begin{lem} \,\ \\
\vspace{-1em}
\begin{enumerate}
\item $\Lambda^3 W_{\Q} \cong [1^3] \oplus [2,1]^{\oplus 2} \oplus [3] \oplus [2,2,-1]^{\oplus 3} \oplus [3,1,-1]$.
\item $\mathcal{L}_F^{\Q}(3) \cong [3] \oplus [2,1]^{\oplus 6} \oplus [1^3] \oplus [2,2,-1]^{\oplus 3} \oplus [3,1,-1]^{\oplus 3}
            \oplus [3,2,-2]^{\oplus 2}$.
\end{enumerate}
\end{lem}
\textit{Proof}.
(1). \,\,From (\ref{eq-ext-dec}),
\[ \Lambda^3 W_{\Q} \cong \Lambda^3 [1] \bigoplus (\Lambda^2 [1,1] \otimes_{\Q} [1,1,-1]) \bigoplus (H_{\Q} \otimes_{\Q} \Lambda^2 [1,1,-1])
    \bigoplus \Lambda^3 [1,1,-1]. \]
Hence, from Lemma {\rmfamily \ref{L-aoi}}, we obtain the required result.

\vspace{0.5em}

(2). \,\,Consider the $\mathrm{GL}(3,\Q)$-equivariant surjective homomorphism $\Lambda^2 W_{\Q} \otimes_{\Q} W_{\Q} \rightarrow \mathcal{L}_F^{\Q}(3)$
whose kernel is $\Lambda^3 W_{\Q}$. On the other hand, we can compute the irreducible decomposition of $\Lambda^2 W_{\Q} \otimes_{\Q} W_{\Q}$
by using Pieri's rule and Lemma {\rmfamily \ref{L-aoi}} as
\[ \Lambda^2 W_{\Q} \otimes_{\Q} W_{\Q} \cong \big{(} [3] \oplus [2,1]^{\oplus 4} \oplus [1^3] \oplus [2,2,-1]^{\oplus 3}
     \oplus [3,1,-1]^{\oplus 2} \oplus [3,2,-2] \big{)}^{\oplus 2}. \]
Thus from Part (1), we obtain the required result.
$\square$

\vspace{0.5em}

Now, we consider the irreducible decomposition of $(R_3/R_4)_{\Q}$.
Consider the $\mathrm{GL}(n,\Q)$-equivariant homomorphism
$[\, , \,]_{\Q} : (R/R_3)_{\Q} \otimes_{\Q} \mathcal{L}_F^{\Q}(1) \rightarrow (R_3/R_4)_{\Q}$
induced from the Lie bracket of $\mathcal{L}_F^{\Q}$.
\begin{pro}\label{P-img}
The map $[\, , \,]_{\Q}$ is injective.
\end{pro}
\noindent
The proof of the above proposition requires patient, lengthy and straightforward computation.
Thus we give the proof in Section {\rmfamily \ref{S-Ap2}} later.

\vspace{0.5em}

In our previous paper \cite{S03}, we have determined the cokernel, and hence the image, of the third rational Johnson homomorphism
$\tau_{3,\Q} : \mathrm{gr}^3(\mathcal{A}_3) \rightarrow H^*_{\Q} \otimes_{\Q} \mathcal{L}_{n,\Q}(4)$. In particular,
we showed that $\mathrm{dim}_{\Q}(\mathrm{gr}_{\Q}^3(\mathcal{A}_3))=43$.
From the fact that $\mathcal{A}_3'(3) = \mathcal{A}_3(3)$ and $\mathcal{A}_3'(4) \subset \mathcal{A}_3(4)$,
we have the surjective homomorphism $\mathrm{gr}_{\Q}^3(\mathcal{A}_3') \rightarrow \mathrm{gr}_{\Q}^3(\mathcal{A}_3)$.
This shows that $\mathrm{dim}_{\Q}(\mathrm{gr}_{\Q}^3(\mathcal{A}_3')) \geq 43$.
Thus, from the exact sequence (\ref{eq-momo}) for $k=3$ and $\mathrm{dim}_{\Q}(\mathcal{L}_F^{\Q}(3))=240$, we see that
$\mathrm{dim}_{\Q}( (R_3/R_4)_{\Q} ) \leq 197$.
On the other hand, by Proposition {\rmfamily \ref{P-img}}, we see that $(R_3/R_4)_{\Q}$ has a subspace
$\mathrm{Im}([\, , \,]_{\Q})$ whose dimension is $162$, and that
\[ \mathrm{dim}_{\Q}((R_3/R_4[R,F])_{\Q}) = \mathrm{dim}_{\Q}(\mathrm{Coker}([\, , \,]_{\Q})) \leq 35. \]

\begin{thm}\label{T-main1}
The module $(R_3/R_4[R,F])_{\Q}$ is an irreducible $\mathrm{GL}(3,\Q)$-representation which is isomorphic to $[3,2,-2]$.
In particular, $\mathrm{dim}_{\Q}((R_3/R_4[R,F])_{\Q}) = 35$.
\end{thm}
\textit{Proof}.
Observe maximal vectors
\[ [K_{312}, \iota_1, K_{312}] = [K_{312}, K_{21}, K_{312}] + [K_{312}, K_{31}, K_{312}], \hspace{1em} [K_{312}, K_{31}, K_{312}] \]
in $\mathcal{L}_F^{\Q}(3)$.
By considering the Hall basis of $\mathcal{L}_F^{\Q}(3)$ with respect to the ordering
\begin{equation}\label{eq-Hall}
 K_{12} < K_{13} < K_{21} < K_{23} < K_{31} < K_{32} < K_{123} < K_{213} < K_{312}
\end{equation}
among the generators of $F$, we see that the above maximal vectors are linearly independent.

\vspace{0.5em}

Since $[K_{312}, \iota_1] \in R$, we see that the maximal vector $[K_{312}, \iota_1, K_{312}]$
belongs to $\mathrm{Im}([\, , \,]_{\Q})$.
On the other hand, by a direct calculation, we see that
\[ r := [K_{312}, K_{31}^{-1}, K_{312}] [ [K_{31}^{-1}, K_{312}], [K_{32}^{-1}, K_{31}^{-1}] ]^{-1} \in R_3, \]
and hence $[K_{312}, K_{31}, K_{312}] = -r \in R_3/R_4$.
Thus $[K_{312}, K_{31}, K_{312}]$ is non-trivial in $(R_3/R_4[R,F])_{\Q}$. Namely, $(R_3/R_4[R,F])_{\Q}$ contains
an irreducible $\mathrm{GL}(3,\Q)$-submodule $[3,2,-2]$ whose dimension is $35$.
Since $\mathrm{dim}_{\Q}((R_3/R_4[R,F])_{\Q}) \leq 35$ as mentioned above, we obtain the required result.
$\square$

\vspace{0.5em}

From Theorem {\rmfamily \ref{T-main1}}, we see that $H^2(\mathrm{IA}_3,\Q)/\mathrm{Im}(\cup_{\Q})$ contains
an irreducible representation $[3,2,-2]^*=[2,-2,-3] \cong D^{-3} \otimes_{\Q} [5,1]$.
On the other hand, as a corollary, we obtain the following.

\begin{cor}\label{C-And}
$\mathcal{A}_3(4)/\mathcal{A}_3'(4)$ is finite.
\end{cor}
\textit{Proof}.
By Theorem {\rmfamily \ref{T-main1}}, we see
\[\begin{split}
   \mathrm{dim}_{\Q}(\mathrm{gr}_{\Q}^3(\mathcal{A}_3')) & = \mathrm{dim}_{\Q}(\mathcal{L}_F^{\Q}(3)) - \mathrm{dim}_{\Q}((R_3/R_4)_{\Q}) \\
   & = 240 -197 = 43.
  \end{split}\]
Namely, the natural surjective homomorphism $\mathrm{gr}_{\Q}^3(\mathcal{A}_3') \rightarrow \mathrm{gr}_{\Q}^3(\mathcal{A}_3)$
is an isomorphism. Hence $\mathcal{A}_3'(4)$ has at most finite index in $\mathcal{A}_3(4)$.
$\square$

\vspace{0.5em}

Finally, we consider the third cohomology group $H^3(\mathrm{IA}_3,\Q)$.
Let $\mathrm{IO}_3$ be the quotient group of $\mathrm{IA}_3$ by the inner automorphism group $\mathrm{Inn}\,F_3$.
Bestvina, Bux and Margalit \cite{Bes} showed that $\mathrm{IO}_3$ has $2$-dimensional Eilenberg-Maclane space, and that
$H_2(\mathrm{IO}_3,\Q)$ is infinitely generated. Hence, by using the spectral sequence of the group extension
\[ 1 \rightarrow \mathrm{Inn}\,F_3 \rightarrow \mathrm{IA}_3 \rightarrow \mathrm{IO}_3 \rightarrow 1, \]
we see that $H^3(\mathrm{IA}_3,\Q)$ is infinitely generated.
The following theorem shows that non-trivial elements in $H^3(\mathrm{IA}_3,\Q)$ cannot be detected by the triple cup product
of the first cohomology group of $\mathrm{IA}_3$.

\begin{thm}\label{T-main2}
The image of the triple cup product
\[ \cup_{\mathrm{IA}_3}^3 : \Lambda^3 H^1(\mathrm{IA}_3, \Q) \rightarrow H^3(\mathrm{IA}_3, \Q) \]
is trivial.
\end{thm}
\textit{Proof}.
Let $\delta : H^0(\mathrm{IA}_n, (R_3/[R,R])_{\Q}^*) \rightarrow H^1(\mathrm{IA}_n, (R/R_3)_{\Q}^*)$
be the connecting homomorphism defined in Subsection {\rmfamily \ref{Ss-coh}}.
By Lemma {\rmfamily \ref{L-conn}}, we see that the restriction of $\delta$ to $(R_3/R_4)_{\Q}^*$ induces the injective homomorphism
\[ \overline{\delta} : (R_3/R_4)_{\Q}^*/(R_3/R_4[F,R])_{\Q}^* \rightarrow H^1(\mathrm{IA}_3, (R/R_3)_{\Q}^*). \]
From Theorem {\rmfamily \ref{T-main1}}, we have
\[\begin{split}
 \mathrm{dim}_{\Q}((R_3 & /R_4)_{\Q}^* /(R_3/R_4[F,R])_{\Q}^*) \\
   & = \mathrm{dim}_{\Q}((R_3/R_4)_{\Q}^*) - \mathrm{dim}_{\Q}((R_3/R_4[F,R])_{\Q}^*) \\
   & = 197- 35 =162.
 \end{split}\]
On the other hand, since
\[ H^1(\mathrm{IA}_3, (R/R_3)_{\Q}^*) \cong H^1(\mathrm{IA}_3, \Q) \otimes_{\Q} (R/R_3)_{\Q}^*, \]
from Lemma {\rmfamily \ref{L-mare}}, we have
\[ \mathrm{dim}_{\Q}(H^1(\mathrm{IA}_3, (R/R_3)_{\Q}^*)) = 162. \]
Therefore, the homomorphism $\overline{\delta}$ is an isomorphism.
Thus $\beta^1 : H^1(\mathrm{IA}_3, (R/R_3)_{\Q}^*) \rightarrow H^3(\mathrm{IA}_3, \Q)$ is zero map.
Then by Proposition {\rmfamily \ref{P-cup3-2}}, we obtain the required result.
$\square$

\section{Proof of Proposition {\rmfamily \ref{P-img}}}\label{S-Ap2}

In this section, we give the proof of Proposition {\rmfamily \ref{P-img}}. The proof is done by patient lengthy hand calculation.
To begin with, we explain how to prove, and show a brief of calculations.

\vspace{0.5em}

By the natural injective homomorphism $R_3/R_4 \hookrightarrow \mathcal{L}_F(3)$, we consider $R_3/R_4$ as a submodule of
$\mathcal{L}_F(3)$.
We show that the homomorphism
\[ [\, , \,] : R/R_3 \otimes_{\Z} \mathcal{L}_F(1) \rightarrow \mathcal{L}_F(3) \]
given by $r \otimes x \mapsto [r,x]$ is injective.
From this, we can obtain Proposition {\rmfamily \ref{P-img}} by tensoring with $\Q$.
In order to do this, we show that the cokernel $\mathcal{L}_F(3)/\mathrm{Im}([\, , \,])$ is a free abelian group of rank
\[ \mathrm{rank}_{\Z}(\mathcal{L}_F(3)) - \mathrm{rank}_{\Z}((R/R_3) \otimes_{\Z} \mathcal{L}_F(1)) = 240 - 162 = 78. \]
From this, we can see that the image $\mathrm{Im}([\, , \,])$ is a free abelian group of rank $162$, and that
the bracket map $[\, , \,]$ is injective.

\vspace{0.5em}

We give basis of $\mathcal{L}_F(1)$ and $R/R_3$ as follows.
First, the set
\[ \mathfrak{E}_1 := \{  K_{12}, K_{13}, K_{21}, K_{23}, K_{31}, K_{32}, K_{123}, K_{213}, K_{312} \} \subset \mathcal{L}_F(1) \]
forms a basis of $\mathcal{L}_F(1)$. Consider the relators

\vspace{0.5em}

\begin{center}
{\renewcommand{\arraystretch}{1.2}
\begin{tabular}{|l|l|} \hline
 {\bf (R1)} & $[K_{ij}, K_{kj}]$ \\ \hline
 {\bf (R2)} & $[K_{ik} K_{jk}, K_{ij}]$ \\ \hline
 {\bf (R3)} & $[K_{ij}K_{kj}, K_{ijk}]$ \\ \hline
 {\bf (R4)} & $[K_{ik}K_{jk}, K_{kij}] (K_{ik}K_{jk}) [K_{ij}, K_{ki}K_{ji}] (K_{ik}K_{jk})^{-1}$ \\ \hline
\end{tabular}}
\end{center}

\vspace{0.5em}

\noindent
in $R$ for distinct $i, j$ and $k$. These elements define a basis of $R/R_3$ given as follows.

\vspace{0.5em}

{\footnotesize
\begin{center}
{\renewcommand{\arraystretch}{1.2}
\begin{tabular}{|l|l|} \hline
 {\bf (R1-1)} & $[K_{32}, K_{12}]$ \\ \hline
 {\bf (R1-2)} & $[K_{23}, K_{13}]$ \\ \hline
 {\bf (R1-3)} & $[K_{31}, K_{21}]$ \\ \hline
 {\bf (R2-1)} & $[K_{23}, K_{12}] + [K_{13}, K_{12}]$ \\ \hline
 {\bf (R2-2)} & $[K_{32}, K_{13}] - [K_{13}, K_{12}]$ \\ \hline
 {\bf (R2-3)} & $[K_{23}, K_{21}] - [K_{21}, K_{13}]$ \\ \hline
 {\bf (R2-4)} & $[K_{32}, K_{31}] - [K_{31}, K_{12}]$ \\ \hline
 {\bf (R2-5)} & $[K_{31}, K_{23}] - [K_{23}, K_{21}]$ \\ \hline
 {\bf (R2-6)} & $[K_{32}, K_{31}] + [K_{32}, K_{21}]$ \\ \hline
 {\bf (R3-1)} & $[K_{123}, K_{32}] + [K_{123}, K_{12}]$ \\ \hline
 {\bf (R3-2)} & $[K_{123}, K_{23}] + [K_{123}, K_{13}]$ \\ \hline
 {\bf (R3-3)} & $[K_{213}, K_{31}] + [K_{213}, K_{21}]$ \\ \hline
 {\bf (R3-4)} & $[K_{213}, K_{23}] + [K_{213}, K_{13}]$ \\ \hline
 {\bf (R3-5)} & $[K_{312}, K_{31}] + [K_{312}, K_{21}]$ \\ \hline
 {\bf (R3-6)} & $[K_{312}, K_{32}] + [K_{312}, K_{12}]$ \\ \hline
 {\bf (R4-1)} & $[K_{312}, K_{23}] + [K_{312}, K_{13}] + [K_{21}, K_{12}] + [K_{31}, K_{12}]$ \\ \hline
 {\bf (R4-2)} & $[K_{123}, K_{31}] + [K_{123}, K_{21}] + [K_{32}, K_{23}] - [K_{23}, K_{12}]$ \\ \hline
 {\bf (R4-3)} & $[K_{213}, K_{12}] + [K_{213}, K_{32}] + [K_{31}, K_{13}] + [K_{21}, K_{13}]$ \\ \hline
\end{tabular}}
\end{center}
}

\vspace{0.5em}

Let $\mathfrak{E}_2 \subset R/R_3$ be the set of the above elements. Then $\mathrm{Im}([\, , \,])$ is generated by 162
elements $[r, x]$ for any $r \in \mathfrak{E}_2$ and $x \in \mathfrak{E}_1$.
Consider again the Hall basis of $\mathcal{L}_F(3)$ consisting of the Magnus generators of $F$ where its order is given by
(\ref{eq-Hall}). The list of the Hall basis is given in Table 1 at the end of this section.
Each element $[K, L, M]$ of the Hall basis of $\mathcal{L}_F(3)$ satisfies $K > L \leq M$.
The elements in Table 1 are arranged in order of $L$. 

\vspace{0.5em}

Now, let $\mathfrak{E}$ be the set of (the coset classes of) the Hall basis in the cokernel $\mathcal{L}_F(3)/\mathrm{Im}([\, , \,])$.
Then $\mathfrak{E}$ forms generators of $\mathcal{L}_F(3)/\mathrm{Im}([\, , \,])$.
For any $r \in \mathfrak{E}_2$ and $x \in \mathfrak{E}_1$, by considering an element $[r, x]$ as a relation
$[r, x] =0$ in $\mathcal{L}_F(3)/\mathrm{Im}([\, , \,])$, we remove some of generators of $\mathcal{L}_F(3)/\mathrm{Im}([\, , \,])$ one by one.
Here we give a few examples.

\vspace{0.5em}

{\bf (1)} $[ {\textbf{(R1-1)}}, K]$.

\vspace{0.5em}

For any $K \in \mathfrak{E}_1$,
we have $[K_{32}, K_{12}, K]=0 \in \mathcal{L}_F(3)/\mathrm{Im}([\, , \,])$, and hence
we remove the generators of type $[K_{32}, K_{12}, K]$ from $\mathfrak{E}$.

\vspace{0.5em}

{\bf (2)} $[ {\textbf{(R1-2)}}, K]$.

\vspace{0.5em}

For any $K \in \mathfrak{E}_1$, we have $[K_{23}, K_{13}, K]=0 \in \mathcal{L}_F(3)/\mathrm{Im}([\, , \,])$.
Namely, if necessary, by using the Jacobi identity, we see
\[ \begin{cases}
     - [K_{13}, K_{12}, K_{23}] + [K_{23}, K_{12}, K_{13}]=0, \hspace{1em} & K = K_{12}, \\
     [K_{23}, K_{13}, K]=0, & K_{13} < K.
   \end{cases}\]
Then we remove the generators of type $[K_{13}, K_{12}, K_{23}]$ and $[K_{23}, K_{13}, K]$ for any $K>K_{13}$ from $\mathfrak{E}$.

\vspace{0.5em}

{\bf (3)} $[ {\textbf{(R1-3)}}, K]$.

\vspace{0.5em}

For any $K \in \mathfrak{E}_1$, by using $[K_{23}, K_{13}, K]=0 \in \mathcal{L}_F(3)/\mathrm{Im}([\, , \,])$.
we remove the generators of type $[K_{21}, K, K_{31}]$ for any $K < K_{21}$, and
$[K_{31}, K_{21}, K]$ for any $K_{21} \leq K$ from $\mathfrak{E}$.

\vspace{0.5em}

{\bf (4)} $[ {\textbf{(R2-1)}}, K]$.

\vspace{0.5em}

For any $K \in \mathfrak{E}_1$, consider $[K_{23}, K_{12}, K] + [K_{13}, K_{12}, K]=0 \in \mathcal{L}_F(3)/\mathrm{Im}([\, , \,])$.
If $K \neq K_{23}$, by using this relation we remove the generators $[K_{13}, K_{12}, K]$ from $\mathfrak{E}$.
If $K = K_{23}$, by the argument in {\bf (2)}, we have
$[K_{13}, K_{12}, K_{23}] = [K_{23}, K_{12}, K_{13}] \in \mathcal{L}_F(3)/\mathrm{Im}([\, , \,])$.
Thus the relation considering here is equivalent to
\[ [K_{23}, K_{12}, K_{23}] + [K_{23}, K_{12}, K_{13}]=0, \]
and hence we remove $[K_{23}, K_{12}, K_{13}]$ from $\mathfrak{E}$.

\vspace{0.5em}

{\bf (5)} $[{\textbf{(R2-2)}}, K]$.

\vspace{0.5em}

For any $K \in \mathfrak{E}_1$, consider $[K_{32}, K_{13}, K] - [K_{13}, K_{12}, K]=0 \in \mathcal{L}_F(3)/\mathrm{Im}([\, , \,])$.
We consider the following four cases.

\vspace{0.5em}

\begin{center}
{\renewcommand{\arraystretch}{1.2}
\begin{tabular}{|l|l|l|l|} \hline
             & $K$         & Relation                   & \\ \hline
 {\bf (i)}   & $K_{12}$    & $[K_{23}, K_{12}, K_{12}] + [K_{23}, K_{12}, K_{32}]=0$  &  {\bf (4)}, {\bf (1)}  \\ \hline
 {\bf (ii)}  & $K_{13}$    & $-[K_{23}, K_{12}, K_{23}] + [K_{32}, K_{13}, K_{13}]=0$ &  {\bf (4)} \\ \hline
 {\bf (iii)} & $K_{23}$    & $[K_{32}, K_{13}, K_{23}] - [K_{32}, K_{13}, K_{13}]=0$  &  {\bf (2)}, {\bf (4)}, {\bf (5)-(ii)}  \\ \hline
 {\bf (iv)}  & otherwise   & $[K_{23}, K_{12}, K] + [K_{32}, K_{13}, K]=0$            &  {\bf (4)} \\ \hline
\end{tabular}}
\end{center}

\vspace{0.5em}

From the above relations, we remove $[K_{23}, K_{12}, K_{12}]$, $[K_{23}, K_{12}, K_{23}]$, $[K_{32}, K_{13}, K_{13}]$ and
$[K_{23}, K_{12}, K]$ for $K \neq K_{12}, K_{13}, K_{23}$ from $\mathfrak{E}$ in turn.

\vspace{0.5em}

Hereafter, we repeat similar computations according to the following order.

\vspace{0.5em}

{\footnotesize
\begin{center}
{\renewcommand{\arraystretch}{1.2}
\begin{tabular}{|l|l|} \hline
                                      & $K$ \\ \hline
 {\bf (6)}  $[{\textbf{(R2-3)}}, K]$  & {\bf (i)} $K=K_{12}$, {\bf (ii)} $K=K_{13}$, {\bf (iii)} $K=K_{31}$, {\bf (iv)} $K \geq K_{21}$ \\ \hline
 {\bf (7)}  $[{\textbf{(R2-4)}}, K]$  & {\bf (i)} $K=K_{12}, K_{13}, K_{21}, K_{23}$, {\bf (ii)} $K \geq K_{31}$ \\ \hline
 {\bf (8)}  $[{\textbf{(R2-5)}}, K]$  & {\bf (i)} $K=K_{31}$, {\bf (ii)} $K \geq K_{23}, \neq K_{31}$, {\bf (iii)} $K = K_{21}, K_{13}, K_{12}$ \\ \hline
 {\bf (9)}  $[{\textbf{(R2-6)}}, K]$  & {\bf (i)} $K \geq K_{31}$, {\bf (ii)} $K = K_{23}, K_{21}, K_{13}, K_{12}$ \\ \hline
 {\bf (10)} $[{\textbf{(R3-1)}}, K]$  & {\bf (i)} $K \geq K_{32}$, {\bf (ii)} $K = K_{31}, K_{23}, K_{21}, K_{13}, K_{12}$ \\ \hline
 {\bf (11)} $[{\textbf{(R3-2)}}, K]$  & {\bf (i)} $K \geq K_{23}$, {\bf (ii)} $K = K_{21}, K_{13}, K_{12}$ \\ \hline
 {\bf (12)} $[{\textbf{(R3-3)}}, K]$  & {\bf (i)} $K \geq K_{31}$, {\bf (ii)} $K = K_{23}, K_{21}, K_{13}, K_{12}$ \\ \hline
 {\bf (13)} $[{\textbf{(R3-4)}}, K]$  & {\bf (i)} $K \geq K_{23}$, {\bf (ii)} $K = K_{21}, K_{13}, K_{12}$ \\ \hline
 {\bf (14)} $[{\textbf{(R3-5)}}, K]$  & {\bf (i)} $K \geq K_{31}$, {\bf (ii)} $K = K_{23}, K_{21}, K_{13}, K_{12}$ \\ \hline
 {\bf (15)} $[{\textbf{(R3-6)}}, K]$  & {\bf (i)} $K \geq K_{32}$, {\bf (ii)} $K = K_{31}, K_{23}, K_{21}, K_{13}, K_{12}$ \\ \hline
 {\bf (16)} $[{\textbf{(R4-1)}}, K]$  & {\bf (i)} $K = K_{23}, K_{31}, K_{32}, K_{123}, K_{213}, K_{312}$,
                                        {\bf (ii)} $K = K_{21}, K_{13}, K_{12}$ \\ \hline
 {\bf (17)} $[{\textbf{(R4-2)}}, K]$  & {\bf (i)} $K = K_{31}, K_{32}, K_{123}, K_{213}, K_{312}$,
                                        {\bf (ii)} $K = K_{23}, K_{21}, K_{13}, K_{12}$ \\ \hline
 {\bf (18)} $[{\textbf{(R4-3)}}, K]$  & {\bf (i)} $K = K_{32}, K_{123}, K_{213}, K_{312}$,
                                        {\bf (ii)} $K = K_{31}, K_{23}, K_{21}, K_{13}, K_{12}$ \\ \hline
\end{tabular}}
\end{center}
}

\vspace{0.5em}

\noindent
In each step, the removed generators are attached the step number in Table 1. For example, $[K_{31}, K_{12}, K_{12}]$ is removed from $\mathfrak{E}$
at the step {\bf (7)-(i)}. Since all of the computations are straightforward and too lengthy to write down here,
we leave it to the high-motivated reader as exercises.

\vspace{0.5em}

Finally, we see that from the above computations, the cokernel $\mathcal{L}_F(3)/\mathrm{Im}([\, , \,])$ is the free abelian group
with basis consisting the elements in Table 1 to which no step number attached.
This completes the proof of Proposition {\rmfamily \ref{P-img}}.

\vspace{0.5em} 

\begin{center}
{\bf Table 1.}
\end{center}

{\tiny
\begin{center}
{\renewcommand{\arraystretch}{1.0}
\begin{tabular}{|l|l|l|l|l|} \hline
   $[K_{13}, K_{12}, K_{12}]$  {\bf (4)}  & $[K_{13}, K_{12}, K_{13}]$  {\bf (4)} & $[K_{13}, K_{12}, K_{21}]$ {\bf (4)}
 & $[K_{13}, K_{12}, K_{23}]$  {\bf (4)}  & $[K_{13}, K_{12}, K_{31}]$  {\bf (4)} \\ \hline
   $[K_{13}, K_{12}, K_{32}]$  {\bf (4)}  & $[K_{13}, K_{12}, K_{123}]$ {\bf (4)} & $[K_{13}, K_{12}, K_{213}]$ {\bf (4)}
 & $[K_{13}, K_{12}, K_{312}]$ {\bf (4)}  & \\ \hline
   $[K_{21}, K_{12}, K_{12}]$             & $[K_{21}, K_{12}, K_{13}]$  {\bf (6i)} & $[K_{21}, K_{12}, K_{21}]$  
 & $[K_{21}, K_{12}, K_{23}]$  {\bf (8iii)}  & $[K_{21}, K_{12}, K_{31}]$  {\bf (3)} \\ \hline
   $[K_{21}, K_{12}, K_{32}]$  {\bf (9ii)}  & $[K_{21}, K_{12}, K_{123}]$           & $[K_{21}, K_{12}, K_{213}]$ {\bf (12ii)}
 & $[K_{21}, K_{12}, K_{312}]$ {\bf (14ii)}  & \\ \hline
   $[K_{23}, K_{12}, K_{12}]$  {\bf (5)}  & $[K_{23}, K_{12}, K_{13}]$  {\bf (4)} & $[K_{23}, K_{12}, K_{21}]$  {\bf (5)}
 & $[K_{23}, K_{12}, K_{23}]$  {\bf (5)}  & $[K_{23}, K_{12}, K_{31}]$  {\bf (5)} \\ \hline
   $[K_{23}, K_{12}, K_{32}]$  {\bf (5)}  & $[K_{23}, K_{12}, K_{123}]$ {\bf (5)} & $[K_{23}, K_{12}, K_{213}]$ {\bf (5)}
 & $[K_{23}, K_{12}, K_{312}]$ {\bf (5)}  & \\ \hline
   $[K_{31}, K_{12}, K_{12}]$  {\bf (7i)}  & $[K_{31}, K_{12}, K_{13}]$  {\bf (7i)} & $[K_{31}, K_{12}, K_{21}]$  {\bf (7i)}
 & $[K_{31}, K_{12}, K_{23}]$  {\bf (7i)}  & $[K_{31}, K_{12}, K_{31}]$  {\bf (7ii)} \\ \hline
   $[K_{31}, K_{12}, K_{32}]$  {\bf (7ii)}  & $[K_{31}, K_{12}, K_{123}]$ {\bf (7ii)} & $[K_{31}, K_{12}, K_{213}]$ {\bf (7ii)}
 & $[K_{31}, K_{12}, K_{312}]$ {\bf (7ii)}  & \\ \hline
   $[K_{32}, K_{12}, K_{12}]$  {\bf (1)}  & $[K_{32}, K_{12}, K_{13}]$  {\bf (1)} & $[K_{32}, K_{12}, K_{21}]$  {\bf (1)}
 & $[K_{32}, K_{12}, K_{23}]$  {\bf (1)}  & $[K_{32}, K_{12}, K_{31}]$  {\bf (1)} \\ \hline
   $[K_{32}, K_{12}, K_{32}]$  {\bf (1)}  & $[K_{32}, K_{12}, K_{123}]$ {\bf (1)} & $[K_{32}, K_{12}, K_{213}]$ {\bf (1)}
 & $[K_{32}, K_{12}, K_{312}]$ {\bf (1)}  & \\ \hline
   $[K_{123}, K_{12}, K_{12}]$  {\bf (10ii)}  & $[K_{123}, K_{12}, K_{13}]$  {\bf (11ii)} & $[K_{123}, K_{12}, K_{21}]$  {\bf (10ii)}
 & $[K_{123}, K_{12}, K_{23}]$              & $[K_{123}, K_{12}, K_{31}]$  \\ \hline
   $[K_{123}, K_{12}, K_{32}]$  {\bf (10i)}  & $[K_{123}, K_{12}, K_{123}]$ {\bf (10i)} & $[K_{123}, K_{12}, K_{213}]$ {\bf (10i)}
 & $[K_{123}, K_{12}, K_{312}]$ {\bf (10i)}  & \\ \hline
   $[K_{213}, K_{12}, K_{12}]$              & $[K_{213}, K_{12}, K_{13}]$  {\bf (13ii)} & $[K_{213}, K_{12}, K_{21}]$   
 & $[K_{213}, K_{12}, K_{23}]$  {\bf (18ii)}  & $[K_{213}, K_{12}, K_{31}]$  {\bf (18ii)} \\ \hline
   $[K_{213}, K_{12}, K_{32}]$  {\bf (18ii)}  & $[K_{213}, K_{12}, K_{123}]$            & $[K_{213}, K_{12}, K_{213}]$  
 & $[K_{213}, K_{12}, K_{312}]$             & \\ \hline
   $[K_{312}, K_{12}, K_{12}]$  {\bf (15ii)}  & $[K_{312}, K_{12}, K_{13}]$  {\bf (16ii)} & $[K_{312}, K_{12}, K_{21}]$   {\bf (15ii)}
 & $[K_{312}, K_{12}, K_{23}]$              & $[K_{312}, K_{12}, K_{31}]$    \\ \hline
   $[K_{312}, K_{12}, K_{32}]$  {\bf (15i)}  & $[K_{312}, K_{12}, K_{123}]$ {\bf (15i)} & $[K_{312}, K_{12}, K_{213}]$  {\bf (15i)}
 & $[K_{312}, K_{12}, K_{312}]$ {\bf (15i)}  & \\ \hline
\end{tabular}

\vspace{1em}

\begin{tabular}{|l|l|l|l|l|} \hline
   $[K_{21}, K_{13}, K_{13}]$  {\bf (6ii)}  & $[K_{21}, K_{13}, K_{21}]$  {\bf (6iv)} & $[K_{21}, K_{13}, K_{23}]$ {\bf (6iv)}
 & $[K_{21}, K_{13}, K_{31}]$  {\bf (3)}  & $[K_{21}, K_{13}, K_{32}]$  {\bf (6iv)} \\ \hline
   $[K_{21}, K_{13}, K_{123}]$ {\bf (6iv)}  & $[K_{21}, K_{13}, K_{213}]$ {\bf (6iv)} & $[K_{21}, K_{13}, K_{312}]$ {\bf (6iv)}
 &  & \\ \hline
   $[K_{23}, K_{13}, K_{13}]$  {\bf (2)}  & $[K_{23}, K_{13}, K_{21}]$  {\bf (2)} & $[K_{23}, K_{13}, K_{23}]$ {\bf (2)}
 & $[K_{23}, K_{13}, K_{31}]$  {\bf (2)}  & $[K_{23}, K_{13}, K_{32}]$  {\bf (2)} \\ \hline
   $[K_{23}, K_{13}, K_{123}]$ {\bf (2)}  & $[K_{23}, K_{13}, K_{213}]$ {\bf (2)} & $[K_{23}, K_{13}, K_{312}]$ {\bf (2)}
 &  & \\ \hline
   $[K_{31}, K_{13}, K_{13}]$             & $[K_{31}, K_{13}, K_{21}]$  {\bf (8i)} & $[K_{31}, K_{13}, K_{23}]$ {\bf (8iii)}
 & $[K_{31}, K_{13}, K_{31}]$             & $[K_{31}, K_{13}, K_{32}]$  {\bf (9ii)} \\ \hline
   $[K_{31}, K_{13}, K_{123}]$            & $[K_{31}, K_{13}, K_{213}]$ {\bf (12ii)} & $[K_{31}, K_{13}, K_{312}]$ {\bf (14ii)}
 &  & \\ \hline
   $[K_{32}, K_{13}, K_{13}]$  {\bf (5)}  & $[K_{32}, K_{13}, K_{21}]$            & $[K_{32}, K_{13}, K_{23}]$ 
 & $[K_{32}, K_{13}, K_{31}]$             & $[K_{32}, K_{13}, K_{32}]$      \\ \hline
   $[K_{32}, K_{13}, K_{123}]$ {\bf (10ii)} & $[K_{32}, K_{13}, K_{213}]$           & $[K_{32}, K_{13}, K_{312}]$ {\bf (15ii)}
 &  & \\ \hline
   $[K_{123}, K_{13}, K_{13}]$              & $[K_{123}, K_{13}, K_{21}]$  {\bf (17ii)} & $[K_{123}, K_{13}, K_{23}]$ {\bf (11i)}
 & $[K_{123}, K_{13}, K_{31}]$  {\bf (11i)}  & $[K_{123}, K_{13}, K_{32}]$  {\bf (11i)} \\ \hline
   $[K_{123}, K_{13}, K_{123}]$ {\bf (11i)}  & $[K_{123}, K_{13}, K_{213}]$ {\bf (11i)} & $[K_{123}, K_{13}, K_{312}]$ {\bf (11i)}
 &  & \\ \hline
   $[K_{213}, K_{13}, K_{13}]$  {\bf (13ii)}  & $[K_{213}, K_{13}, K_{21}]$  {\bf (13ii)} & $[K_{213}, K_{13}, K_{23}]$ {\bf (13i)}
 & $[K_{213}, K_{13}, K_{31}]$  {\bf (13i)}  & $[K_{213}, K_{13}, K_{32}]$  {\bf (13i)} \\ \hline
   $[K_{213}, K_{13}, K_{123}]$ {\bf (13i)}  & $[K_{213}, K_{13}, K_{213}]$ {\bf (13i)} & $[K_{213}, K_{13}, K_{312}]$ {\bf (13i)}
 &  & \\ \hline
   $[K_{312}, K_{13}, K_{13}]$  {\bf (16ii)}  & $[K_{312}, K_{13}, K_{21}]$   & $[K_{312}, K_{13}, K_{23}]$
 & $[K_{312}, K_{13}, K_{31}]$              & $[K_{312}, K_{13}, K_{32}]$    \\ \hline
   $[K_{312}, K_{13}, K_{123}]$             & $[K_{312}, K_{13}, K_{213}]$  & $[K_{312}, K_{13}, K_{312}]$ 
 &  & \\ \hline
\end{tabular}

\vspace{1em}

\begin{tabular}{|l|l|l|l|l|} \hline
   $[K_{23}, K_{21}, K_{21}]$  {\bf (8iii)}  & $[K_{23}, K_{21}, K_{23}]$  {\bf (8ii)} & $[K_{23}, K_{21}, K_{31}]$ {\bf (6iii)}
 & $[K_{23}, K_{21}, K_{32}]$  {\bf (8ii)}  & $[K_{23}, K_{21}, K_{123}]$ {\bf (8ii)} \\ \hline
   $[K_{23}, K_{21}, K_{213}]$ {\bf (8ii)}  & $[K_{23}, K_{21}, K_{312}]$ {\bf (8ii)} & & & \\ \hline
   $[K_{31}, K_{21}, K_{21}]$  {\bf (3)}  & $[K_{31}, K_{21}, K_{23}]$  {\bf (3)} & $[K_{31}, K_{21}, K_{31}]$ {\bf (3)}
 & $[K_{31}, K_{21}, K_{32}]$  {\bf (3)}  & $[K_{31}, K_{21}, K_{123}]$ {\bf (3)} \\ \hline
   $[K_{31}, K_{21}, K_{213}]$ {\bf (3)}  & $[K_{31}, K_{21}, K_{312}]$ {\bf (3)} & & & \\ \hline
   $[K_{32}, K_{21}, K_{21}]$  {\bf (9ii)}  & $[K_{32}, K_{21}, K_{23}]$  {\bf (9ii)} & $[K_{32}, K_{21}, K_{31}]$ {\bf (9i)}
 & $[K_{32}, K_{21}, K_{32}]$  {\bf (9i)}  & $[K_{32}, K_{21}, K_{123}]$ {\bf (9i)} \\ \hline
   $[K_{32}, K_{21}, K_{213}]$ {\bf (9i)}  & $[K_{32}, K_{21}, K_{312}]$ {\bf (9i)} & & & \\ \hline
   $[K_{123}, K_{21}, K_{21}]$             & $[K_{123}, K_{21}, K_{23}]$  {\bf (17ii)} & $[K_{123}, K_{21}, K_{31}]$ {\bf (17ii)}
 & $[K_{123}, K_{21}, K_{32}]$  {\bf (17ii)} & $[K_{123}, K_{21}, K_{123}]$  \\ \hline
   $[K_{123}, K_{21}, K_{213}]$            & $[K_{123}, K_{21}, K_{312}]$  & & & \\ \hline
   $[K_{213}, K_{21}, K_{21}]$  {\bf (12ii)} & $[K_{213}, K_{21}, K_{23}]$   & $[K_{213}, K_{21}, K_{31}]$
 & $[K_{213}, K_{21}, K_{32}]$  {\bf (18ii)} & $[K_{213}, K_{21}, K_{123}]$  \\ \hline
   $[K_{213}, K_{21}, K_{213}]$            & $[K_{213}, K_{21}, K_{312}]$  & & & \\ \hline
   $[K_{312}, K_{21}, K_{21}]$  {\bf (14ii)} & $[K_{312}, K_{21}, K_{23}]$             & $[K_{312}, K_{21}, K_{31}]$  {\bf (14i)}
 & $[K_{312}, K_{21}, K_{32}]$  {\bf (14i)} & $[K_{312}, K_{21}, K_{123}]$ {\bf (14i)} \\ \hline
   $[K_{312}, K_{21}, K_{213}]$ {\bf (14i)} & $[K_{312}, K_{21}, K_{312}]$ {\bf (14i)} & & & \\ \hline
\end{tabular}

\vspace{1em}

\begin{tabular}{|l|l|l|l|l|} \hline
   $[K_{31}, K_{23}, K_{23}]$             & $[K_{31}, K_{23}, K_{31}]$            & $[K_{31}, K_{23}, K_{32}]$
 & $[K_{31}, K_{23}, K_{123}]$ {\bf (11ii)} & $[K_{31}, K_{23}, K_{213}]$ {\bf (12ii)} \\ \hline
   $[K_{31}, K_{23}, K_{312}]$ {\bf (14ii)}  & & & & \\ \hline
   $[K_{32}, K_{23}, K_{23}]$             & $[K_{32}, K_{23}, K_{31}]$            & $[K_{32}, K_{23}, K_{32}]$  
 & $[K_{32}, K_{23}, K_{123}]$ {\bf (10ii)} & $[K_{32}, K_{23}, K_{213}]$ {\bf (18ii)} \\ \hline
   $[K_{32}, K_{23}, K_{312}]$ {\bf (15ii)}  & & & & \\ \hline
   $[K_{123}, K_{23}, K_{23}]$  {\bf (11ii)} & $[K_{123}, K_{23}, K_{31}]$            & $[K_{123}, K_{23}, K_{32}]$  
 & $[K_{123}, K_{23}, K_{123}]$            & $[K_{123}, K_{23}, K_{213}]$ \\ \hline
   $[K_{123}, K_{23}, K_{312}]$  & & & & \\ \hline
   $[K_{213}, K_{23}, K_{23}]$   & $[K_{213}, K_{23}, K_{31}]$            & $[K_{213}, K_{23}, K_{32}]$
 & $[K_{213}, K_{23}, K_{123}]$  & $[K_{213}, K_{23}, K_{213}]$ \\ \hline
   $[K_{213}, K_{23}, K_{312}]$  & & & & \\ \hline
   $[K_{312}, K_{23}, K_{23}]$  {\bf (16i)} & $[K_{312}, K_{23}, K_{31}]$  {\bf (16i)}   & $[K_{312}, K_{23}, K_{32}]$  {\bf (16i)}
 & $[K_{312}, K_{23}, K_{123}]$ {\bf (16i)} & $[K_{312}, K_{23}, K_{213}]$ {\bf (16i)}   \\ \hline
   $[K_{312}, K_{23}, K_{312}]$ {\bf (16i)} & & & & \\ \hline
\end{tabular}
}
\end{center}}

\newpage

{\tiny
\begin{center}
{\renewcommand{\arraystretch}{1.0}
\begin{tabular}{|l|l|l|l|l|} \hline
   $[K_{32}, K_{31}, K_{31}]$  {\bf (16ii)} & $[K_{32}, K_{31}, K_{32}]$            & $[K_{32}, K_{31}, K_{123}]$  {\bf (10ii)}
 & $[K_{32}, K_{31}, K_{213}]$            & $[K_{32}, K_{31}, K_{312}]$ {\bf (15ii)} \\ \hline
   $[K_{123}, K_{31}, K_{31}]$  {\bf (17i)} & $[K_{123}, K_{31}, K_{32}]$  {\bf (17i)} & $[K_{123}, K_{31}, K_{123}]$  {\bf (17i)}
 & $[K_{123}, K_{31}, K_{213}]$ {\bf (17i)} & $[K_{123}, K_{31}, K_{312}]$ {\bf (17i)} \\ \hline
   $[K_{213}, K_{31}, K_{31}]$  {\bf (12i)} & $[K_{213}, K_{31}, K_{32}]$  {\bf (12i)} & $[K_{213}, K_{31}, K_{123}]$  {\bf (12i)}
 & $[K_{213}, K_{31}, K_{213}]$ {\bf (12i)} & $[K_{213}, K_{31}, K_{312}]$ {\bf (12i)} \\ \hline
   $[K_{312}, K_{31}, K_{31}]$    & $[K_{312}, K_{31}, K_{32}]$    & $[K_{312}, K_{31}, K_{123}]$
 & $[K_{312}, K_{31}, K_{213}]$   & $[K_{312}, K_{31}, K_{312}]$   \\ \hline
\end{tabular}

\vspace{1em}

\begin{tabular}{|l|l|l|l|} \hline
   $[K_{123}, K_{32}, K_{32}]$    & $[K_{123}, K_{32}, K_{123}]$    & $[K_{123}, K_{32}, K_{213}]$   
 & $[K_{123}, K_{32}, K_{312}]$   \\ \hline
   $[K_{213}, K_{32}, K_{32}]$  {\bf (18i)}  & $[K_{213}, K_{32}, K_{123}]$  {\bf (18i)}  & $[K_{213}, K_{32}, K_{213}]$   {\bf (18i)}
 & $[K_{213}, K_{32}, K_{312}]$ {\bf (18i)}  \\ \hline
   $[K_{312}, K_{32}, K_{32}]$    & $[K_{312}, K_{32}, K_{123}]$    & $[K_{312}, K_{32}, K_{213}]$
 & $[K_{312}, K_{32}, K_{312}]$   \\ \hline
\end{tabular}

\vspace{1em}

\begin{tabular}{|l|l|l|} \hline
   $[K_{213}, K_{123}, K_{123}]$  & $[K_{213}, K_{123}, K_{213}]$  & $[K_{213}, K_{123}, K_{312}]$ \\ \hline
   $[K_{312}, K_{123}, K_{123}]$  & $[K_{312}, K_{123}, K_{213}]$  & $[K_{312}, K_{123}, K_{312}]$ \\ \hline
   $[K_{312}, K_{213}, K_{213}]$  & $[K_{312}, K_{213}, K_{312}]$  &  \\ \hline

\end{tabular}

}
\end{center}}

\section{Acknowledgments}\label{S-Ack}

The author would like to express his sincere gratitude to Professor Shigeyuki Morita for discussing this research many times
after his lecture for undergraduate students at Tokyo University of Science.
His insatiable thirst for mathematics and lots of indisputable leading works have been the intense longing and the perfect model of the author
as a mathematician.
The author has been inspired by his mathematical spirits, and supported by his unfailing encouragements for more than ten years.

\vspace{0.5em}

The author would like to thank Andrew Putman for giving me usuful comments about the Andreadakis-Johnson filtration of $\mathrm{Aut}\,F_n$.
The author also would like to thank the referee for the careful reading and helpful comments for
the original paper.

\vspace{0.5em}

This work is supported by JSPS KAKENHI Grant Number 16K05155.

\end{document}